\documentclass{article}

\usepackage{amsthm}
\usepackage{amsmath}
\usepackage{amssymb}
\usepackage{amsfonts}
\usepackage{authordate1-4}

\def\N{\mathbb{N}}
\def\P{\mathbb{P}}
\def\R{\mathbb{R}}
\def\S{\mathbb{S}}
\def\E{\mathbb{E}}

\def\Z{\mathbb{Z}}

\def\Me{\mathbb{\mathcal{M}}}

\def\C{\mathcal{C}}
\def\Pr{\mathcal{P}}

\def\B{\mathcal{B}}
\def\ED{\mathcal{E}}

\def\F{\mathcal{F}}
\def\D{\mathcal{D}}
\def\I{\mathcal{I}}
\def\T{\mathbb{T}}
\def\<{\big\langle}
\def\>{\big\rangle}
\def\Osc{\operatorname{Osc}}

\def\eref#1{(\ref{#1})}

\newtheorem{Theorem}{Theorem}[section]

\newtheorem{Lemma}[Theorem]{Lemma}
\newtheorem{Corollary}[Theorem]{Corollary}

\newtheorem{Proposition}[Theorem]{Proposition}
\newtheorem{Conjecture}[Theorem]{Conjecture}

\theoremstyle{remark}
\newtheorem{Remark}[Theorem]{Remark}
\theoremstyle{definition}

\begin{document}

\title{Anomalous slow diffusion from perpetual homogenization.\protect\footnotetext{Received May 20, 2001; revised }\protect\footnotetext{AMS 1991 {\it{Subject Classification}}. Primary 60J60; secondary  , 35B27, 34E13, 60G44, 60F05, 31C05.} \protect\footnotetext{{\it{Key words and phrases}}. Multi scale homogenization, anomalous diffusion, diffusion on fractal media, heat kernel,  subharmonic, exponential martingale inequality, Davies's conjecture, periodic operator.}}         
\author{Houman Owhadi \protect\footnote{LATP, UMR CNRS 6632, CMI, Universit\'{e} de Provence , owhadi@cmi.univ-mrs.fr}}        

\maketitle

\abstract{This paper is concerned with the asymptotic behavior of solutions of  stochastic differential equations
$dy_t=d\omega_t -\nabla V(y_t) dt$, $y_0=0$. When $d=1$ and  $V$ is not periodic but obtained as a superposition  of an infinite number of
 periodic potentials with geometrically increasing periods ($V(x) = \sum_{k=0}^\infty U_k(x/R_k)$, where $U_k$ are smooth functions of period 1, $U_k(0)=0$, and
$R_k$ grows exponentially fast with $k$) we can show that $y_t$ has an anomalous slow behavior
 and we obtain quantitative estimates on the anomaly using and developing the tools of homogenization. Pointwise estimates are based on a new analytical inequality for sub-harmonic functions. When $d\geq 1$ and $V$ is periodic,  quantitative estimates are obtained on the heat kernel of $y_t$, showing the rate at which homogenization takes place. The latter result proves Davies's conjecture
 and is based on a quantitative estimate for the Laplace transform of martingales
that can be used to obtain similar results for periodic elliptic generators}

\section{Introduction}
It is now well known that natural Brownian Motions on various disordered or complex structures are anomalously slow.\\
These mechanisms of the slow diffusion for instance are well understood for very regular strictly self-similar fractals. The archetypical specific example of a deep problem being the one solved in \cite{BB97} on the Sierpinski Carpet (which is infinitely ramified, a codeword for hard to understand rigorously: for a survey on diffusions on fractals we refer to \cite{Ba98}, for an alternative approach  to \cite{Os95} and for the random Sierpinski Carpet  to \cite{HaKu98}). It appears that the main feature is the existence of an infinite number of scales of obstacle (with proper size) for the diffusion.\\
It is our object to show that one can implement the common idea that this last feature (infinitely many scales) is the key to the possibility of anomalous diffusion, in a general context using the tools of homogenization.\\
The strategy of the proof might appear paradoxical: it is not a priori very sensible to try to prove that the diffusion is anomalous by the use of homogenization theory which is a vast mathematical machine destined to prove an opposite result, i.e a central limit theorem and thus  normal diffusion. But it will be shown that when the homogenization process  is not finished, an anomalous behavior whose characteristics are controlled by homogenization theory  might appear.\\
This paper will focus on the sub-diffusive behavior in dimension one (subsection \ref{MRsdbeha1MR1}), which will allow
the introduction of a concept of differentiation between spatial scales that can be applied to a more general framework.\\
The proof of the anomaly of the exit times is based on a new quantitative analytical inequality for sub-harmonic functions (subsection \ref{MRsdbeha1MR3}) that is linked with stability properties
of elliptic divergence form operators.\\
The extension of those results to higher dimensions has been done in \cite{BeOw00b} and to the super-diffusive case in \cite{BeOw00c} and \cite{Ow00b}.\\
The control of the anomalous
heat kernel tail is based on sharp quantitative estimates for the Laplace transform of a martingale. These estimates  allow us to put into evidence the rate at which homogenization takes place on the behavior of the heat kernel of an elliptic generator in any dimension (subsection \ref{MRsdbeha1MR2}). The quantitative control of the heat kernel in homogenization theory outside any asymptotic regime has been recognized as difficult and important  \cite{Nor97}. For instance, this problem is at the center of Davies' conjecture emphasized as "well beyond existing results"  \cite{Dav93}. With theorem \ref{SMSpCocorkjgptperme1} we give a proof of that conjecture in any dimension for elliptic operators with only bounded coefficients.

\subsection{History}
The idea of associating homogenization  (or renormalization) on large number of scales with the anomaly of a physical system has already been applied from an heuristic point of view to several physical models.\\
Maybe one of the oldest of such applications is to Differential Effective Medium theories which was first proposed by Bruggeman to calculate the conductivity of a two-component composite structure formed by successive substitutions (\cite{Bru35} and \cite{AIP77}) and generalized in \cite{Nor85} to materials with more than two phases. For instance this theory has been applied to compute the anomalous electrical and acoustic properties of fluid-saturated sedimentary rocks \cite{SeScCo81}. More recently this problem has been analyzed from a rigorous point of view in \cite{Av87} and  \cite{Ko95};  in \cite{AlBr96} and  \cite{JiKo99}.\\
The heuristic application of this idea to prove the anomalous behavior of  diffusion seems to have been done only for the super-diffusive case that is to say for a diffusion evolving among a large number of divergence-free drifts. Maybe this is explained by the strong motivation to explore convective transports in turbulent flows which are known to be characterized by a large number of scales of eddies. The first observation was empirical: in  \cite{Ric26}  Richardson empirically conjectured that the diffusion coefficient $D_\lambda$ in turbulent air depends on the scale length $\lambda$ of the measurement.
More recently physicists and mathematicians have started to investigate on the super-diffusive phenomenon (from both heuristic and rigorous points of view)  using the tools of homogenization or renormalization (the first cousin of multi-scale homogenization): we refer to      \cite{AvMa90}; \cite{GlZh92}, \cite{Ave96};   \cite{Bha99}; \cite{FaPa94}; \cite{FaKo01}.\\

\subsection{The model}
Let us consider in dimension one a Brownian motion with a drift given by the gradient of a potential $V$, i.e. the
solution of the stochastic differential equation:
\begin{equation}\label{IntModelsubdiffstochdiffequ}
  dy_t = d\omega_t - \nabla V(y_t) dt , \quad y_0=0.
\end{equation}
The multi-scale potential  $V$ is given by a sum of infinitely many periodic functions with
(geometrically) increasing periods:
\begin{equation}\label{Modsubfracuinfty}
V=\sum_{n=0}^\infty U_n(\frac{x}{R_n})
\end{equation}

In this formula we have two important ingredients: the potentials $U_k$
and the scale parameters $R_k$. We will now describe the hypothesis we
make on these two items of our model.
\begin{enumerate}
\item \underline{Hypotheses on the potentials $U_k$}\\
We will assume that
\begin{equation}\label{jhshdddikuou1}
U_k \in C^\infty(\T)
\end{equation}
\begin{equation}\label{jhshdikuou1}
U_k(0)=0
\end{equation}
Here $C^\infty(\T^d)$ denotes the space of smooth functions
on the torus $\T:=\R/\Z$.
We will also assume that the first derivate of the $U_k$ are uniformly
bounded, i.e.
\begin{equation}\label{ModsubContUngradUn}
K_1:=\sup_{k\in \N} \sup_{x\not=y}|U_k(x)-U_k(y)|/|x-y|.
<\infty
\end{equation}
We will also need the notation
\begin{equation}\label{ModsubContUngradUnGer}
K_0:=\sup_{k\in \N} \Osc(U_k)
\end{equation}
where the oscillation of $U_k$ is given by $\Osc(U):=\sup U - \inf
U$.\\
We write $D(U_k)$ for the effective diffusivities associated to
the  potentials $U_k$: if $z_t$ is the solution of
$dz_t=d\omega_t-\nabla U_k (z_t) dt$ it is well known \cite{Ol94}
that as $\epsilon \downarrow 0$, $\epsilon z_{t/\epsilon^2}$
converges in law towards a Brownian Motion with covariance matrix
$D(U_k)$ given by
\begin{equation}
D(U_n)=\big(\int_{\T}e^{2U_n(x)}dx \int_{\T}e^{-2U_n(x)}dx\big)^{-1}.
\end{equation}
We also assume that the effective diffusivity matrices of the $U_k$'s
are uniformly bounded away from $0$ and $1$.
\begin{equation}\label{ModsubDiffCondUniDUn}
\lambda_{\min}=\inf_{n\in \N} D(U_n)>0 \quad \text{and}\quad  \lambda_{\max}=\sup_{n\in \N} D(U_n)<1.
\end{equation}
\item \underline{Hypotheses on the scale parameters $R_k$}\\
$R_k$ is a spatial scale parameter growing exponentially fast with
$k$, more precisely we will assume that $R_0=r_0=1$ and that the
ratios between scales defined by (we write $\N^{*}$ the set of
integers different from $0$)
\begin{equation}\label{jahgsvagvjh6761}
 r_k= R_k/R_{k-1}\in \N^{*}
\end{equation}
for $k\geq 1$, are integers uniformly bounded away from $1$ and
$\infty$: we will denote by
\begin{equation}\label{Modsubboundrnrhonmin}
 \rho_{\min}:=\inf_{k\in \N^*} r_k \quad\text{and}\quad
\rho_{\max}:=\sup_{k\in \N^*} r_k
\end{equation}
and assume that
\begin{equation}\label{Modsubboundrnrhonmin2}
 \rho_{\min} \geq 2\quad\text{and}\quad \rho_{\max}  < \infty.
\end{equation}
\end{enumerate}

Since $\|\nabla V\|_\infty <\infty$ it is well known that the solution of \eref{IntModelsubdiffstochdiffequ} exists; is unique up to sets of measure $0$ with respect to the Wiener measure and is a strong Markov continuous Feller process.
\begin{Remark}
Note that if   $\forall n,U_n\in\{W_1,\ldots,W_p\}$,  the $(W_i)$ being non constant, then the conditions \eref{ModsubDiffCondUniDUn} and \eref{ModsubContUngradUn} are trivially satisfied.
\end{Remark}

\section{Main results}
\subsection{Sub-diffusive behavior}\label{MRsdbeha1MR1}
Our first objective is to show that the solution of \eref{IntModelsubdiffstochdiffequ} is abnormally slow and the asymptotic sub-diffusivity will be characterized in three ways:
\begin{itemize}
\item as an anomalous behavior of the expectation of $\tau(0,r)$ (the
  exit time  from a ball of radius $r$, for $r\to\infty$, i.e.
$\mathbb E_0[\tau(0,r)] \sim  r^{2+\nu}$).
\item as an anomalous behavior of the variance at time $t$,
  i.e.  $\mathbb E_0[y_t^2] \sim t^{1-\nu}$ as $t\to\infty$.
\item  as an anomalous (non-Gaussian) behavior of the tail of the
  transition probability of the process.
\end{itemize}
More precisely  there exists a constant $\rho_0(K_0,K_1,\lambda_{\max})$ such that
\begin{Theorem}\label{IntSMOnedSuHTcoa6}
If $\rho_{\min}>\rho_0$ and $\tau(0,r)$ is the exit time associated to the solution of \eref{IntModelsubdiffstochdiffequ} then there exists a constant $C_1$ depending on $K_0,K_1$ such that
\begin{equation}\label{estimtau}
\E_0[\tau(0,r)]=r^{2+\nu_1(r)+\epsilon(r)}
\end{equation}
where $\epsilon(r)\rightarrow 0$ as $r\rightarrow \infty$ and
\begin{equation}\label{IntSMOnedSuHTcoa6eq1}
0<-\frac{\ln \lambda_{\max}}{\ln \rho_{\max}}-\frac{C_1}{\rho_{\min}\ln \rho_{\max}}
 \leq \nu_1(r)\leq -\frac{\ln \lambda_{\min}}{\ln \rho_{\min}}+\frac{C_1}{\rho_{\min}\ln
 \rho_{\min}}.
\end{equation}
\end{Theorem}

\begin{Theorem}\label{IntSMOnedsunmsqdthhvh873}
If $\rho_{\min}>\rho_0$ and $y_t$ is a solution of \eref{IntModelsubdiffstochdiffequ} then for there exists a constant $C_2$ depending on $K_0,K_1$ and a time $t_0$ depending on $K_1,\rho_{\min},\rho_{\max},\lambda_{\max}$ such that for $t>t_{0}$
\begin{equation}\label{estimcarre}
\E[y_t^2]=t^{1-\frac{\nu_2(t)}{2}}
\end{equation}
where
\begin{equation}\label{IntSMOnedsunmsqdthhvh873eq1}
0<- \frac{\ln \lambda_{\max}}{\ln \rho_{\max}}- \frac{C_{2}}{\ln
\rho_{\min} \ln \rho_{\max}} \leq \nu_2(t)\leq - \frac{\ln
\lambda_{\min}}{\ln \rho_{\min}}+ \frac{C_{2}}{(\ln
\rho_{\min})^2}.
\end{equation}
\end{Theorem}

\begin{Theorem}\label{IntSMOnetrprdeborathi81}
If $\rho_{\min}>\rho_0$ and $y_t$ is a solution of \eref{IntModelsubdiffstochdiffequ} then there exist constants $C_5$ depending on $K_0,K_1,R_2$, $C_3$ on $K_0,K_1,\rho_{\min}$, $C_4,C_6,C_7$ on $K_0,K_1$ such that if $t,h>0$ and
\begin{equation}\label{IntSMOnedsutrprdethbishe1}
\frac{t}{h} \geq C_{5} \quad \text{and} \quad \frac{h^2}{t}\geq C_{3}(\frac{t}{h})^{\frac{\ln \lambda_{\max}}{2\ln \rho_{\max}}+\frac{C_{4}}{(\ln \rho_{\min})^2}}
\end{equation}
then
\begin{equation}\label{IntSMOnedsutrprdethbishe3}
\ln \P[|y_t|\geq h] \leq -C_{6}\frac{h^2}{t}(\frac{t}{h})^{\nu_3}
\end{equation}
with
\begin{equation}\label{IntSMOnedsutrprdethbishe4}
\nu_3=-\frac{\ln \lambda_{\max}}{\ln \rho_{\max}}-\frac{C_{7}}{\ln
\rho_{\min} \ln \rho_{\max}}>0.
\end{equation}
\end{Theorem}
\begin{Remark}
The second condition in \eref{IntSMOnedsutrprdethbishe1} is really needed since the leading exponent associated to $(t/h)$ is  ($\frac{\ln \lambda_{\max}}{2\ln \rho_{\max}}$), i.e. half the one  associated to $\nu_3$. This condition corresponds to a frontier with a heat kernel diagonal regime.
\end{Remark}

\subsubsection{Description of the proofs}
 Before discussing the results further we want to describe the proof. A perpetual homogenization process takes place over the infinite number of scales $0,\ldots,n,\ldots$. The idea is to distinguish, when one tries to estimate \eref{estimtau}, \eref{estimcarre} or \eref{IntSMOnedsutrprdethbishe3}, the smaller scales  which have already been homogenized ($0,\ldots,n_{ef}$ called effective scales), the bigger scales  which have not had a visible influence on the diffusion ($n_{dri},\ldots,\infty$ called drift scales because they will be replaced by a constant drift in the proof) and some intermediate scales that manifest their particular shapes in the behavior of the diffusion ($n_{ef}+1,\ldots,n_{dri}-1=n_{ef}+n_{per}$ called perturbation scales because they will enter in the proof as a perturbation of the homogenization process over the smaller scales). To estimate \eref{estimtau} for instance, if
one considers the periodic approximation of the potential
\begin{equation}
  V_0^n(x) = \sum_{k=0}^n U_k(x/R_k)
\end{equation}
the corresponding process $y_t^{(n)}$  will have an asymptotic
(\emph{homogenized}) variance \cite{Ol94}
\begin{equation}\label{sdkjcbkjbckjbrkjfbbfffs1}
D(V_0^n)=\big(\int_{\T}e^{2V_0^n(R_n x)}dx \int_{\T}e^{-2V_0^n(R_n x)}dx\big)^{-1}
\end{equation}
 $D(U_0)$ is  smaller than $1$ and because of the geometric growth of the periods $R_n$ and a minimal separation between them (i.e. $\rho_{\min}>\rho_0$), $D(V_0^n)$ decreases exponentially fast in $n$. \\
By homogenization theory, $y_t^n$ is characterized by a mixing length $\xi_m(V_0^n)\sim R_n$ such that if one writes $\tau^n$ its associated exit times then for $r>\xi_m(V_0^n)$
\begin{equation}
  \E_0[\tau^n(0,r)] \sim \frac{r^2}{D(V_0^{n})}.
\end{equation}
Writing $n_{ef}(r)=\sup\{n: R_n\le r\}$ one proves that $\E_0[\tau(0,r)]\sim \E_0[\tau^{n_{ef}(r)}(0,r)]$ by showing the stability of $\E_0[\tau(0,r)]$ under the influence of $V_{n_{ef}(r)+1}^\infty=\sum_{k=n_{ef}(r)+1}^\infty U_k(x/R_k)$. This control is based on a new analytical inequality which
 shall be described in the sequel and allows to obtain that
\begin{equation}\label{Newanlaineqjn1constau}
\begin{split}
\E_0[\tau^{n_{ef}(r)}(0,r)]& e^{-6 \Osc_r(V_{n_{ef}(r)+1}^\infty)} \\&\leq \E_0[\tau(0,r)]\leq
  \E_0[\tau^{n_{ef}(r)}(0,r)] e^{6
  \Osc_r(V_{n_{ef}(r)+1}^\infty)}.
\end{split}
\end{equation}
In these inequalities $\Osc_r(V_{n_{ef}(r)+1}^\infty)$ stands for $\sup_{B(0,r)}V_{n_{ef}(r)+1}^\infty-\inf_{B(0,r)}V_{n_{ef}(r)+1}^\infty$ and is controlled by $$\Osc_r(V_{n_{ef}(r)+1}^\infty)\leq \Osc(U_{n_{ef}(r)+1}) + \|\nabla V_{n_{ef}(r)+2}^\infty\|_\infty r$$ i.e. $n_{ef}(r)+1$ acts as a perturbation scale and $n_{ef}(r)+2,\ldots,\infty$ as drift scales.
 From this
\begin{equation}
  \E_0[\tau(0,r)] \sim \frac{r^2}{D(V_0^{n_{ef}(r)})}.
\end{equation}
 Thus, if
\begin{equation*}
 - \lim \inf_{r\to\infty} \frac{1}{\ln r} \ln D(V_0^{n_{ef}(r)}) > 0
\end{equation*}
one has sub-diffusivity, in the sense as defined above.\\
The proof of \eref{estimcarre} follows similar lines by the
introduction mixing times $\tau_m(V_0^n)$ and visibility times
$\tau_v(V_p^\infty)$ (such that for
$\tau_m(V_0^n)<t<\tau_v(V_p^\infty)$, $V_p^\infty$ has not a real
influence on the behavior of the diffusion $y_t$ and $V_0^n$ has
been homogenized). Then choosing $n_{ef}(t)=\sup\{n:
\tau_m(V_0^n)\leq t\}$ one obtains that
\begin{Proposition}\label{propcompcaryt2ujhb1}
Let $\nu_2(t)$ the function associated to \eref{estimcarre}, one
has for $t>t_{K_1,\rho_{\min},\rho_{\max},\lambda_{\max}}$
\begin{equation}\label{eqmorethcoto}
\nu_{ef}(t) (1-\frac{C_{K_1}}{\ln \rho_{\min}}) \leq \nu_2(t)\leq \nu_{ef}(t) (1+\frac{C_{K_1}}{\ln \rho_{\min}})
\end{equation}
\begin{equation}
\nu_{ef}(t)=\frac{\ln \frac{1}{\lambda_{ef}(t)}}{\ln \rho_{ef}(t)}\quad \text{with}
 \quad \rho_{ef}^{n_{ef}}=R_{n_{ef}} \text{and} \quad
 \lambda_{ef}^{n_{ef}+1}=D(V_0^{n_{ef}}).
\end{equation}
\end{Proposition}
This proposition shows than this separation between scales is more than a conceptual tool, it does reflect the underlying phenomenon. Indeed  the anomalous function $\nu_2(t)$ is given in the first order in $1/(\ln \rho_{\min})$ by the number of effective scales by $\E[y_t^2]\sim t D(V_0^{n_{ef}})$, and in this approximation $\nu_2(t)\sim \nu_{ef}(t)$ where $\nu_{ef}(t)$ corresponds to a medium in which the ratios $r_n$ and the effectives diffusivities $D(U_n)$ have been replaced by  their geometric mean over the $n_{ef}+1$ effective scales. The origin of the constant $C_{K_1}/(\ln \rho_{\min})$ in \eref{eqmorethcoto} is the perturbation scales. More precisely, one has to fix the drift scales by $n_{dri}(t)=\inf\{n: \tau_v(V_n^\infty)\geq t\}$, and in general there is a gap between $n_{ef(t)}$ and $n_{dri}(t)$, the scales $U_n$ situated in this gap manifest their particular shape in the behavior of $\nu_2(t)$ and since no hypothesis have been made on those shapes one has to take into account their influence as a perturbation.

\paragraph{}
One may notice that in many papers on diffusions on fractals (see e.g. \cite{Ba98} section 3) obtaining estimates on hitting times is essentially the key to the whole problem and the same is true here: this  strategy has been adapted in \cite{BeOw00b}. In this paper we have chosen to not use this strategy in order to put an emphasis on the role played by the never-ending homogenization process taking place on these diffusions on fractals. Indeed one might wonder why the estimates of the behavior of  Brownian Motions on fractals are of the form
\begin{equation}\label{SMOranbetsfeq1}
\E[y_t^2] \sim t^\frac{2}{d_w},
\end{equation}
\begin{equation}\label{SMOranbetsfeq2}
\E[\tau(0,r)] \sim r^{d_w},
\end{equation}
\begin{equation}\label{SMOranbetsfeq3}
\ln p(t,x,y) \sim -
\big(\frac{|x-y|^{d_w}}{t}\big)^\frac{1}{d_w-1}.
\end{equation}
One explanation  is given here by the number of effective scales hidden in  the estimates \eref{SMOranbetsfeq1},  \eref{SMOranbetsfeq2} and \eref{SMOranbetsfeq3}. Let us assume the model to be self similar (for all $k$, $r_k=\rho$ and $U_k=U$, $D(U_k)=\lambda$). In the table below we have summarized formulae giving (in the first approximation in $1/\ln \rho$) the number of effective scales and the formulae linking them with those anomalous estimates (appearing in the proof,   the influence of the perturbation scales will be neglected). This  gives three values of $d_w$ corresponding to  \eref{SMOranbetsfeq1}, \eref{SMOranbetsfeq2}, \eref{SMOranbetsfeq3} and the interesting point is to compare them.

\begin{table}
\begin{center}
\begin{tabular}{l| c c c}
 & $\E_0[y_t^2]$ & $\E_0[\tau(0,y)]$ & $\ln \P_0[y_t\geq h]$\\ \\\hline \\
$n_{ef}$ & $\frac{\ln t}{2 \ln \rho}$ & $\frac{\ln r}{\ln \rho}$ & $\frac{\ln \frac{t}{h}}{ \ln \frac{\rho}{\lambda^\frac{1}{2}}}$ \\ \\
Heuristic & $t \lambda^{n_{ef}}$ & $\frac{r^2}{\lambda^{n_{ef}}}$ & $-\frac{h^2}{t\lambda^{n_{ef}}}$\\ \\
Anomaly & $t^{\frac{2}{d_{w,1}}}$ & $r^{d_{w,2}}$ & $-\big(\frac{h^{d_{w,3}}}{t}\big)^{\frac{1}{d_{w,3}-1}} $\\ \\
$d_{w,i}$ & $\frac{2}{1+\frac{\ln \lambda}{2 \ln \rho}}$ & $2-\frac{\ln \lambda}{\ln \rho}$ & $1+\frac{1}{1+\frac{\ln \lambda}{\ln \rho-\frac{1}{2}\ln \lambda}}$
\end{tabular}
\end{center}
\end{table}

Let us observe that the multi-scale homogenization techniques gives back the right forms for the mean squared displacement, the exit times and the transition probability densities;  they are explained by the number of  scales  which homogenization can be considered as complete associated to each observation. Moreover $d_{w,1}, d_{w,2}$ and $d_{w,3}$ are  equal up the first order approximation in $1/\ln \rho$ nevertheless they are not equal and this is not surprising. Indeed when $\rho$ is small the second order term in $1/(\ln \rho)^2$ can not be neglected since the perturbation scales becomes more and more dominant (and  the influence of the perturbation scales is of the order of $1/(\ln \rho)^2$).

\subsubsection{Strong overlap between the spatial scales}
The anomaly is based on a minimal separation between spatial scales i.e. $\rho_{\min}>\rho_0$ and one might wonder what happens below this boundary. The answer will be given on a self similar case, i.e. $V$ is said to be self similar if for all $n$, $U_n=U$ and $\rho_{\min}=\rho_{\max}=\rho$.
\begin{Theorem}\label{IntSMOnedSuHTcoa5}
If the potential $V$ in \eref{IntModelsubdiffstochdiffequ} is self similar. Then for all $\rho\geq 2$
\begin{equation}
\E_0[\tau(0,r)]=r^{2+\nu(r)}
\end{equation}
with
\begin{equation}\label{IntSUMoeqRenukhz81}
\nu(r)=\frac{\Pr_{\rho}(2U)+\Pr_{\rho}(-2U)}{\ln \rho}+\epsilon(r)
\end{equation}
with $\epsilon(r)\rightarrow 0$ as $r\rightarrow \infty$.
\end{Theorem}
Here $\Pr_{\rho}$ is the topological pressure associated to the shift operator $s_{\rho}: x\in \T \rightarrow \rho x \in \T$ (see \eref{ProToVaFoPrenent} for its definition).\\
Using the convexity properties of the topological pressure one has $\Pr_{\rho}(2U)+\Pr_{\rho}(-2U)\geq 0$ and \begin{Proposition}\label{IntSMOnedSuHTcoahhhbj5}
$\Pr_{\rho}(2U)+\Pr_{\rho}(-2U)=0$
if and only if
\begin{equation}
\lim_{n\rightarrow \infty}
\Big\|\frac{1}{n}\sum_{k=0}^{n-1}\big(U(\rho^k
x)-\int_{\T^d}U(x)dx\big)\Big\|_{\infty}=0.
\end{equation}
\end{Proposition}
From this one deduces that for the simple example $U(x)=\sin(x)-\sin(81x)$,
$\E[\tau(0,r)]$ is anomalous (sub-diffusive $\sim r^{2+\nu}$ with $\nu>0$)
for $\rho \in \{2\}\cup \{4,\ldots,26\}\cup\{28,\ldots 80\}\cup \{82,\ldots,+\infty \}$
and normal ($\sim r^2$) for $\rho=3,27,81$.\\
Thus if $U$ is not a constant function, there exists $\rho_0(K_0,K_1,D(U))$ such that
for $\rho>\rho_0$, $y_t$  has a clear anomalous behavior ($\E_0[\tau(0,r)]\sim r^{2+\nu}$ with $\nu>0$)
 but in the interval $(1,\rho_0]$
both cases are possible: $y_t$ may show a normal  or an anomalous behavior according
to the value of the ratio between
scales $\rho$ and the regions of normal behavior (characterized by proposition
\ref{IntSMOnedSuHTcoahhhbj5}) might be separated by regions of anomalous behavior. \\
What creates this phenomenon is a strong overlap or interaction between scales: that is
why the region $(1,\rho_0)$ will be called "overlapping ratios", i.e. in this region the
fluctuation of $V$ at a size $\xi>0$ is not represented by a single $U_n(x/R_n)$ but by
 several ones and to characterize the behavior of $y_t$ in that region one must
introduce additional parameters describing the shapes of the fluctuations $U_n$, elsewhere a normal or a sub-diffusive behavior are both possible.

\subsection{Davies's conjecture and quantitative estimates on rate of convergence towards the limit process in homogenization}\label{MRsdbeha1MR2}
The proof of theorem \ref{IntSMOnetrprdeborathi81} has not been described yet.
 The strategy is still to distinguish effective, perturbation and drift scales
nevertheless it is not obvious to determine how many scales have been homogenized in the estimation of $\P_0(y_t \geq h)$.
The answer is directly linked with the rate at which the transition probability densities associated
with a periodic elliptic operator do pass from a large deviation behavior to a \emph{homogenized} behavior.\\
Consider for instance in any dimension $d\geq 1$, $U\in L^\infty(\T^d)$ and the Dirichlet form
\begin{equation}\label{eqdirforfeeduf1}
\ED(f,f)=1/2\int_{R^d}|\nabla f(x)|^2
\frac{e^{-2U(x)}}{\int_{\T^d}e^{-2U(z)}dz}dx,\quad
f\in\D[\ED]=H^1(\R^d).
\end{equation}
Write $p(t,x,y)$ its associated heat kernel with respect to
\begin{equation}
m_U(dx)= \frac{e^{-2U(x)}}{\int_{\T^d}e^{-2U(z)}dz}dx
\end{equation}
 the invariant measure associated to \eref{eqdirforfeeduf1}. Note that when $U$ is smooth the associated operator can be written
$L=1/2\Delta - \nabla U \nabla$ and it is well known that
\begin{itemize}
\item \textbf{Large deviation regime}: for $|x-y|>>t$ the paths of the diffusion concentrate on the geodesics and
\begin{equation}
\ln p(t,x,y)\sim  -\frac{|x-y|^2}{2t}.
\end{equation}
\item \textbf{Heat kernel diagonal regime}: for  $|x-y|^2<<t$, the behavior is fixed by the diagonal of the heat kernel and
\begin{equation}
p(t,x,y)\sim \frac{C_0(x)}{t^\frac{d}{2}}.
\end{equation}
\end{itemize}
Davies conjectured that (we refer to \cite{Dav93}, he considers periodic operators of divergence form nevertheless the idea remains unchanged) that $p(t,x,y)$ should have a \emph{homogenized} behavior ($\ln p(t,x,y)\sim -(x-y)D(U)^{-1}(x-y)/(2t)$) \emph{for $t$ large enough}.\\
J. R. Norris  \cite{Nor97} has shown that the homogenized behavior of the heat kernel $p(t,x,y)$ corresponding to a periodic operator on the torus $\T^d$ (dimension $d$ side $1$) starts at least for $t \ln t >> |x-y|^2$ (with  $|x-y|^2<<t$ ); in this paper it will be shown that it starts for $t >> |x-y|$ in any dimension.\\
This allows to complete the picture describing the behavior of $p(t,x,y)$
\begin{itemize}
\item \textbf{Homogenization regime}: for $1<<|x-y|<<t$ and $|x-y|^2>>t$, homogenization takes place and
\begin{equation}
\ln p(t,x,y)\sim -|x-y|^2_{D^{-1}(U)}/(2t)
\end{equation}
with
\begin{equation}
|x-y|^2_{D^{-1}(U)}:={^t(x-y)}D(U)^{-1}(x-y).
\end{equation}
\end{itemize}
More precisely we will prove that

\begin{Theorem}\label{SMSpCocorkjgptperme1}
Consider $p(t,x,y)$ the heat kernel associated to the Dirichlet form
 \eref{eqdirforfeeduf1} with respect to the measure $m_U$. Then there exist  constants $C,C_2$ depending only on $d$ and $\Osc(U)$ such that for
\begin{equation}
C |x-y| < t,\quad C \sqrt{t}< |x-y|,\quad C<|x-y|
\end{equation}
one has
\begin{equation}\label{dsjddddjsh99111d1}
p(t,x,y)\leq \frac{1}{(2\pi t)^\frac{d}{2} \big(\operatorname{\det(D(U))}\big)^\frac{1}{2}}\exp\big(-(1-E)|y-x|^2_{D^{-1}(U)}/(2t)\big)
\end{equation}
\begin{equation}\label{dsjddddjsh99111d2}
p(t,x,y)\geq \frac{1}{(2\pi t)^\frac{d}{2} \big(\operatorname{\det(D(U))}\big)^\frac{1}{2}}
 \exp\big(-(1+E)|y-x|^2_{D^{-1}(U)}/(2t)\big).
\end{equation}
With
\begin{equation}
E(t,x,y):=  C_2
\big(\frac{|x-y|}{t}+\frac{\sqrt{t}}{|x-y|}\big)\leq \frac{1}{10}.
\end{equation}
\end{Theorem}
Theorem \ref{SMSpCocorkjgptperme1} proves  Davies's conjecture, moreover $E(t,x,y)$ acts as a quantitative error term putting into evidence the rate at which homogenization takes place for the heat kernel, and  it  also acts as the inverse of a distance from  the domains associated to  the large deviation regime and the heat kernel diagonal regime. Observe that all the constants do depend only on $d$ and $\Osc(U)$. It is straightforward to extend those estimates to any periodic elliptic operator. They can be liked to results obtained by A. Dembo \cite{Dem96} for discrete martingales with bounded jumps based on moderate deviations techniques.

\subsubsection{A note on the proof of theorem \ref{IntSMOnetrprdeborathi81}}
Those estimates basically say that the homogenized behavior of the heat kernel associated to a periodic medium of period $R$ starts for $t> R |x-y|$. Thus in the proof of theorem \ref{IntSMOnetrprdeborathi81} the number of the smaller scales that can be considered as homogenized is fixed by $n_{ef}(t/h)=\sup_{n}\{R_n \leq t/h\}$, which (assume $D(U_n)=\lambda$ and $R_n=\rho^n$ for simplification) leads to an anomaly of the form
\begin{equation}\label{eqart1sugtapsierptxy1}
\ln \P(y_t\geq h)\leq -C \frac{h^2}{t \lambda^{n_{ef}(t/h)}}\sim -C \frac{h^2}{t} (\frac{t}{h})^{-\frac{\ln \lambda}{\ln \rho}}\sim - C \big(\frac{|x-y|^{d_w}}{t}\big)^\frac{1}{d_w-1}
\end{equation}
with $d_w\sim 2- \frac{\ln \lambda}{\ln \rho}$. The equation \eref{eqart1sugtapsierptxy1} suggests that the origin of the anomalous shape of the heat kernel for the reflected Brownian Motion on the Sierpinski carpet can be explained by the formula linking the number of effective scales and the ratio $t/h$.\\
The first condition in \eref{IntSMOnedsutrprdethbishe1} can be translated into "homogenization has started on at least the first scale" and the second one into "the heat kernel associated to \eref{IntModelsubdiffstochdiffequ} is far from its diagonal regime"  (one can have $h^2/t<<1$ before reaching that regime, this is explained by the slow down of the diffusion).

\subsubsection{A quantitative inequality for exponential martingales}
The core of the proof of  theorems \ref{IntSMOnetrprdeborathi81} and \ref{SMSpCocorkjgptperme1} is an inequality giving a quantitative estimates for the Laplace transform of a martingale:\\
Consider $M_t$  a continuous square integrable $\F_t$ adapted martingale such that $M_0=0$ and for $\lambda \in \R$, $\E[e^{\lambda M_t}]<\infty$. Assume that there exists a function $f: \R^+ \rightarrow \R^+$ such that for all $t_2 > t_1 \geq 0$ one has a.s.
\begin{equation}\label{eqcontrbracthemartj8bv1}
 \E[\int \limits_{t_1}^{t_2}d<M,M>_s | \F_{t_1}] \leq \int \limits_{0}^{t_2-t_1}
 f(s)ds
\end{equation}
with $f(s) =f_1$ for $s<t_0$ and $f(s)=f_2$ for $s \geq t_0$ with $t_0 > 0$ and $0 < f_2 < f_1$.\\
\begin{Theorem}\label{SMPrtoshcolatrho}
Let $M_t$ be the martingale described above.
\begin{enumerate}
\item  for  all $0<|\lambda|<\big(2e(f_1-f_2)t_0 \big)^{-\frac{1}{2}}$ one has
\begin{equation}\label{eqthmartpreqjhk8kj1}
\E [\exp(\lambda M_t)]  \leq  e^{3(1-1/g(\lambda))}\exp(\frac{g(\lambda)}{2} \lambda^2 f_2 t)
\end{equation}
with $g(\lambda)=\frac{1}{1-\lambda^2(f_1-f_2)t_0 e}$ that verifies $1\leq g \leq 2$
\item for all $0<\nu<\big(2e(f_1-f_2)t_0\big)^{-1}$ one has
\begin{equation}\label{eqthmartpreqjhk8kj2}
\E [\exp(\nu <M,M>_t)] \leq \exp(\nu f_2 t) \frac{\exp\big(\nu
t_0(f_1-f_2) \big)}{((f_1-f_2) \nu t_0)^2}.
\end{equation}
\end{enumerate}
\end{Theorem}
This theorem uses the knowledge on the conditional behavior of the quadratic variation of a martingale to upper bound its Laplace transform, and it is well known that a quantitative control on the Laplace transform leads to a quantitative control on the heat kernel tail.
 The condition $\lambda$ \emph{small enough}  marks the boundary between the large deviation regime and the homogenization regime. A direct application of the key theorem is the following result.
\begin{Corollary}\label{SMPrtoaphoupgoma}
Let $M_t$ be the martingale given in theorem \ref{SMPrtoshcolatrho}.\\
Write $C_1=\big(2e(f_1-f_2)t_0 \big)^\frac{1}{2}/f_2$. For $r=\frac{C_1 x}{ t} < 1$
one has
\begin{equation}
\P(M_t \geq x) \leq
e^{\frac{3}{2}r^2}\exp\big(-(1-r^2)\frac{x^2}{2f_2t} \big).
\end{equation}
\end{Corollary}
This corollary gives a quantitative control on the tail of the law of $M_t$ from the asymptotic behavior of its conditional brackets.

\subsection{An analytical inequality for sub-harmonic functions}\label{MRsdbeha1MR3}
The stability property \eref{Newanlaineqjn1constau} is based on the following analytical inequality:
\begin{Theorem}\label{IntthmSMAldTiger4}
Let $\Omega$ be an open bounded subset of $\R$ ($d=1$),
for $\lambda \in C^\infty(\overline{\Omega})$ such that $\lambda>0$ on $\overline{\Omega}$ and $\phi,\psi \in C^2(\overline{\Omega})$ null on $\partial \Omega$ and both sub harmonic with respect to the operator $-\nabla(\lambda \nabla)$, one has
\begin{equation}\label{equuhbusdciujjh1}
\int_{\Omega}\lambda(x) |\nabla \phi(x).\nabla \psi(x)|\,dx \leq 3
\int_{\Omega}\lambda(x) \nabla \phi(x).\nabla \psi(x)\,dx.
\end{equation}
\end{Theorem}
The constant $3$ in this theorem is the optimal one. We believe that this inequality might also be true in higher dimensions, i.e.:
\begin{Conjecture}\label{IntSMAldTiger4}
For $\Omega \subset \R^d$an open subset with smooth boundary, there exist a constant $C_{d,\Omega}$ depending only on the dimension of the space and the open set  such that
for $\lambda \in C^\infty(\overline{\Omega})$ such that $\lambda>0$ on $\overline{\Omega}$ and $\phi,\psi \in C^2(\overline{\Omega})$ null on $\partial \Omega$ and both sub harmonic with respect to the operator $-\nabla(\lambda \nabla)$, one has
\begin{equation}
\int_{\Omega}\lambda(x) |\nabla \phi(x).\nabla \psi(x)|\,dx \leq
C_{d,\Omega} \int_{\Omega}\lambda(x) \nabla \phi(x).\nabla
\psi(x)\,dx.
\end{equation}
\end{Conjecture}
This conjecture  is  equivalent to the stability of the Green functions of divergence form elliptic operators under a deformation. More precisely write $G_{\lambda}$ the Green function associated to $-\nabla (\lambda \nabla)$ with Dirichlet conditions on $\partial \Omega$.
\begin{Proposition}\label{Propartpreqgrfost}
The conjecture \ref{IntSMAldTiger4} is true with the constant $C_{d,\Omega}$ if and only if for all $\lambda,\mu$ bounded and strictly positive on $\overline{\Omega}$
\begin{equation}\label{Newanghgsinh897hj81}
\big(\sup_{\overline{\Omega}}
\max(\frac{\mu}{\lambda},\frac{\lambda}{\mu})\big)^{-C_{d,\Omega}}
\leq \frac{G_{\mu}(x,y)}{G_{\lambda}(x,y)} \leq
\big(\sup_{\overline{\Omega}}
\max(\frac{\mu}{\lambda},\frac{\lambda}{\mu})\big)^{C_{d,\Omega}}.
\end{equation}
\end{Proposition}
\begin{Remark}
Thus it would be interesting to prove it since it would allow to obtain sharp quantitative estimates on
the comparison of elliptic operators with non Laplacian principal part. By proposition \ref{Propartpreqgrfost} it is easy to check that conjecture \ref{IntSMAldTiger4} implies Harnack inequality. One might think that one would be able to obtain \eref{Newanghgsinh897hj81} using Aronson's estimates and keeping track of the dependence of the constants in the Harnack inequality, but this is not the case since Harnack inequality is an isotropic inequality and \ref{Newanghgsinh897hj81} compares in an optimal way Green functions of operators which can be strongly anisotropic.\\
 Let us remind that the Harnack inequality associated to the operator $L=-\nabla \lambda \nabla$ says that for all $L$-harmonic functions $u$ in $B(0,r)$ one has $$\sup_{x\in B(0,r/2)} u(x)\leq C_L \inf_{x\in B(0,r/2)},$$ where the optimal constant $C_L$ grows towards infinity as $\sup \lambda/\inf \lambda \rightarrow \infty$ whereas the constant associated to conjecture \ref{IntSMAldTiger4} is independent of $\lambda$.
That is why the Harnack inequality strategy, which has already been used to obtain quantitative results for the comparison with the Laplace operator (we refer to \cite{St1}, \cite{Anc97}, \cite{GrWi82} and \cite{Pin92}) allows to obtain
\begin{equation}\label{Newanghgsinh897hjsss81}
\frac{G_{\lambda}(x,y)}{G_{0}(x,y)}\leq C_H
\end{equation}
but with a constant $C_H$ exploding like $C_d\exp\big(C_d(\sup\lambda /\inf \lambda)^{C_d}\big)$
\end{Remark}
\begin{Remark}
Since the conjecture
is true in dimension one with $C_{d,\Omega}=3$ (this constant is an homotopy invariant), it is through proposition \ref{Propartpreqgrfost} that one obtains stability property \eref{Newanlaineqjn1constau}.
\end{Remark}
\begin{Remark}
It is easy to deduce from Theorem \ref{IntthmSMAldTiger4} that if $\Omega$ is a bounded open subset of $\R^d$ and
 $\phi, \psi$ are both convex or both concave functions on $\Omega$ and null on $\partial \Omega$, then
\begin{equation}
\int_{\Omega}|\nabla_x \phi(x).\nabla_x \psi(x)| \,dx \leq 3
\int_{\Omega}\nabla_x \phi(x).\nabla_x \psi(x) \,dx.
\end{equation}
\end{Remark}
\begin{Remark}
The conjecture \ref{IntSMAldTiger4} (theorem \ref{IntthmSMAldTiger4} when $d=1$) has an interesting signification (and consequences) in the framework of electrostatic theory, we refer  to the
chapter 13 of \cite{Owh01}.
\end{Remark}

\subsection{Remark: fast separation between scales}
The feature that distinguishes a strong slow behavior from a weak one is the rate at which spatial scales do separate. Indeed one can follow the proofs given above, changing the condition $\rho_{\max}<\infty$ into  $R_n=R_{n-1}[\rho^{n^\alpha}/R_{n-1}]$ ($\rho,\alpha>1$)  and $\lambda_{\max}=\lambda_{\min}=\lambda<1$ to obtain
\begin{itemize}
\item A weak slow behavior of the exit times
\begin{equation}
C_1 r^{2}e^{g(r)} \leq \E_0[\tau(0,r)] \leq C_2 r^{2}e^{g(r)}
\end{equation}
with $g(r)= (\ln r)^{\frac{1}{\alpha}}(\ln 1/\lambda) (\ln \rho)^{-\frac{1}{\alpha}}$.
\item A weak slow behavior of the mean squared displacement
\begin{equation}
C_1 t e^{-f(t)} \leq \E_0[y_t^2] \leq C_2 t e^{-f(t)}
\end{equation}
with $f(t)= (\ln t)^{\frac{1}{\alpha}}(\ln 1/\lambda)(2\ln \rho)^{-\frac{1}{\alpha}}(1+\epsilon(t))$
\item A weak slow behavior of the heat kernel tail: for $h>0$, $C_1 <t/h<C_2 h$
\begin{equation}\label{IntSMOnedsutrprdethbishe3b}
\P[y_t\geq h] \leq C_{3} e^{-C_{4}\frac{h^2}{t}k(\frac{t}{h})}
\end{equation}
with $k(x)=\lambda^{-(\frac{x}{\ln \rho})^\frac{1}{\alpha}(1+\epsilon(x))}$
\end{itemize}
And as $\alpha \downarrow 1$ the behavior of the solution of \eref{IntModelsubdiffstochdiffequ} pass from weakly anomalous to strongly anomalous.

\section{Proofs}
\subsection{Davies's conjecture and quantitative estimates on
rate of convergence towards the limit process in homogenization}
\subsubsection{Quantitative control of the Laplace transform of a martingale: theorem \ref{SMPrtoshcolatrho}}
The core of the proof of the anomalous heat kernel tail (theorem \ref{IntSMOnetrprdeborathi81}) and the
quantitative estimates on the heat kernel associated to an elliptic operator (theorem \ref{SMSpCocorkjgptperme1}) is  theorem \ref{SMPrtoshcolatrho} that will be proven
in this sub subsection.\\
Let $M_t$ be the martingale described in theorem \ref{SMPrtoshcolatrho}.
Let $q>1$, using H\"{o}lder inequality and Ito formula it is easy to obtain that with $h_q=\frac{q^2}{2(q-1)}$
\begin{equation}\label{eqprbrdemhqjh8ksp30}
\E [\exp(\lambda M_t)] \leq  \E [\exp(h_q\lambda^2 <M,M>_t)]^{\frac{1}{q}}
\end{equation}
Thus the quantitative control of the Laplace transform of the martingale shall follow from this control on its bracket.\\
Write  $\mu=\frac{t}{t_0}$ ($[\mu]$ shall stand for the integer
part of $\mu$). Using H\"{o}lder inequality and  the control
\eref{eqcontrbracthemartj8bv1} one obtains for $1<z<\infty$
\begin{equation}\label{eqprbrxxdemhqjh8ksp31}
\begin{split}
\E [\exp(h_q\lambda^2 <M,M>_t)]^{\frac{1}{q}} \leq&  \E [\exp(z
h_q \lambda^2 <M,M>_{[\mu]t_0})]^\frac{1}{z q}
  \\& \exp\big( (h_q/q)\lambda^2 (t-[\mu]t_0)f_1 \big).
\end{split}
\end{equation}
Then by taking the limit $z\downarrow 1$, one easily obtains that
\begin{equation}\label{eqprbrdemhqjh8ksp31}
\begin{split}
\E [\exp(h_q\lambda^2 <M,M>_t)]^{\frac{1}{q}} \leq&  \E [\exp(h_q \lambda^2 <M,M>_{[\mu]t_0})]^\frac{1}{q}
  \\& \exp\big((h_q/q)\lambda^2 (t-[\mu]t_0)f_1 \big).
\end{split}
\end{equation}
Write $a=\frac{f_2}{f_1}$, we will need the following lemma
\begin{Lemma}\label{dsjhdhdbdhbhh888111}
Let $M_t$ be the martingale described in theorem \ref{SMPrtoshcolatrho} and $\eta >0$,  for $a=f_2/f_1$ and $\mu=t/t_0$ one has
\begin{equation}\label{secprnbreqlksvikd72n80}
\E [\exp(\eta <M,M>_t)]
  \leq 1+ \sum \limits_{n=1}^{+\infty} \frac {(\eta f_1 t_0)^n}{n!} \sum
\limits_{0 \leq m \leq n \wedge \mu}(\mu -m)^n C^{m}_n (a-1)^m.
\end{equation}
\end{Lemma}
\begin{proof}
By the Taylor expansion of the exponential one obtains
\begin{equation}\label{secprnbreqlksvikd72n81}
\exp(\eta <M,M>_t) = 1+ \sum \limits_{n=1}^{+\infty}
\eta^n W_n
\end{equation}
with
 $W_n=\int 1(0<t_1< \cdots <t_n<t) d<M,M>_{t_1} \cdots d<M,M>_{t_n}$. Using the control \eref{eqcontrbracthemartj8bv1}
on the conditional brackets of the martingale it is easy to obtain by
induction on the integrand and the Markov property that
\[ \E[W_n]
\leq \int\limits_{u_i>0} 1(0<u_1 + \cdots +u_n<t) f(u_1) \cdots f(u_n) du_1 \cdots du_n \]
Combining this with \eref{secprnbreqlksvikd72n81} and using the fact that $f(s) \leq f_1 g(s/t_0)$
 with $g(z)=1(z<1)+a 1(z \geq 1)$ one obtains that
\begin{equation}\label{secprnbreqlksvikd72n82}
\E[\exp(\eta <M,M>_t)] \leq 1+ \sum \limits_{n=1}^{+\infty} (\eta
f_1 t_0)^n G_n
\end{equation}
with  $G_n=\int\limits_{z_i>0} 1(0<z_1 + \cdots +z_n<\mu)
 \prod \limits_{k=1}^{n}  (1(z_k<1)+a 1(z_k \geq 1)) dz_1 \cdots dz_n$.
Developing the product in $G_n$ one obtains by integration, induction and straightforward combinatorial
computation that
\[\begin{split}
G_n=\frac{1}{n!} \sum_{0 \leq m \leq \mu \wedge n} C^m_n (\mu-m)^n
(a-1)^m.
\end{split}\]
Which leads to \eref{secprnbreqlksvikd72n80} by the inequality \eref{secprnbreqlksvikd72n82}.

\end{proof}
Using lemma \ref{dsjhdhdbdhbhh888111} one obtains
\[\begin{split}\E [&\exp(h_q\lambda^2 <M,M>_{[\mu]t_0})] \\&\leq  \sum \limits_{n=0}^{+\infty}
 \frac {(h_q \lambda^2 f_1 t_0)^n}{n!} \sum \limits_{0 \leq m \leq n \wedge [\mu]}([\mu] -m)^n C^{m}_n (a-1)^m \end{split}\]
Changing the order of summation, one obtains
\begin{equation}\label{eqconbespprjhlk71}
\begin{split}
\E [\exp(h_q\lambda^2& <M,M>_{[\mu]t_0})] \leq \exp(h_q \lambda^2
f_1 t_0 [\mu]) \\&\sum \limits_{ 0 \leq m \leq [\mu]}  \frac
{([\mu] -m)^m (h_q (a-1)\lambda^2 f_1 t_0)^m}{m!}.
\end{split}
\end{equation}
Now we will need the following lemma
\begin{Lemma}
for $-\frac{1}{e}<y<0$
\begin{equation}\label{eqconbespprjhlk72}
\sum \limits_{ 0 \leq m \leq [\mu]}  \frac {([\mu] -m)^m y^m}{m!}
\leq  \frac{\exp(y [\mu])}{y^2}.
\end{equation}
\end{Lemma}

\begin{proof}
Put $-\frac{1}{e}<x<0$ and write for $n \in \N$,
\[I_n=\sum_{0 \leq m \leq n} \frac{x^m}{m!}(n-m)^m\]
It will be shown here that $\forall p\in \N^*, \; \forall n\in\N$
\begin{equation}\label{eqrefzobimartinegbh87h1}
I_n \leq  \big(u_p(x) \big)^{-n} \big(1-u_p\exp(xu_p)\big)^{-1}
\end{equation}
where $u_p$ the increasing sequence defined by $u_0=0$ and
$u_{p+1}=\exp(-x u_p)$ and converging to $y_0$ the smallest
positive solution of $y \exp(xy)=1$.\\
The inequality \eref{eqconbespprjhlk72}
is then obtained for $u_p(y)=u_2(y)=\exp(-y)$ and using $\exp(-y)-1\geq -y$ and $-\frac{1}{e}<y<0$.\\
Write $y_1=\inf \big\{y>0 \; :\; y\exp(|x|y)=1 \big\}$ (note that $0<y_1<1$) and consider for $-y_1<y<y_1$ the function
\[f\,:\,y\rightarrow \big(1-y \exp(xy)\big)^{-1}\]
By Taylor expansion, for $y\in (-y_1,y_1)$, $f(y)=\sum_{n=0}^{+\infty}y^n \sum_{m=0}^{+\infty}\frac{(nxy)^m}{m!}$
and since\\ $\sum_{0\leq n,m \leq +\infty} y^n \frac{(n|x|y)^m}{m!}=\frac{1}{1-y\exp(|x|y)}<\infty$
with a normal convergence of the series, the order of the limits can be changed, which leads to
\[\begin{split}
f(y)&=\sum_{m=0}^{+\infty}\frac{(nxy)^m}{m!} \sum_{n=0}^{+\infty} n^m y^n
=\sum_{m=0}^{+\infty}\frac{x^m}{m!} \sum_{n=m}^{+\infty} (n-m)^m y^n
\\
&=\sum_{n=0}^{+\infty} y^n \sum_{m=0}^{n} (n-m)^m
\frac{x^m}{m!}=\sum_{n=0}^{+\infty} y^n I_n.
\end{split}\]
It follows that $\forall n\in \N, \; I_n=\frac{f^{(n)}(0)}{n!}$
Now, for $-\frac{1}{e}<x<0$; the constant $y_0=\inf \big\{y>0 \;:\; y\exp(xy)=1\big\}$ does exist and
 $\forall y \in ]-y_1,y_0[,\;\forall n ,\;f^{(n)}(y)\geq 0$ (thus $I_n\geq 0$).\\
Thus from the classical theorem of Taylor expansion,
 the series \\$\sum_{n=0}^{+\infty} y^n \frac{f^{(n)}(0)}{n!}$ converges towards f for $y\in]-y_1,y_0[$ and in that interval
\[\sum_{n=0}^{\infty}y^n I_n=\big(1-y\exp(xy)\big)^{-1}\]
From which one deduces that  $\forall y\in ]0,y_0[\; \forall n\in \N \;
I_n\leq y^{-n} \big(1-y\exp(xy)\big)^{-1}$\\
On the other hand if one considers the sequence $u_0=0,\; u_{p+1}=\exp(-x u_p)$
 then it is an exercise to show that $u_p$ is increasing and will converge towards $y_0$,
which leads to \eref{eqrefzobimartinegbh87h1}.
\end{proof}

Applying \eref{eqconbespprjhlk72}  to \eref{eqconbespprjhlk71} with $y=h_q(a-1) \lambda^2 f_1 t_0$ one obtains that for
$ 0<|\lambda|
<\big(eh_q(f_1-f_2)t_0 \big)^{-\frac{1}{2}}$
one has
\begin{equation}\label{eqconbespprjhlk73}
\E [\exp(h_q\lambda^2 <M,M>_{[\mu]t_0})]^\frac{1}{q} \leq
\exp(\frac{h_q}{q} \lambda^2 f_2 t_0 [\mu])\,(h_q(1-a) \lambda^2
f_1 t_0)^{-\frac{2}{q}}.
\end{equation}
Writing $\nu=\lambda^2 h_q$ and combining \eref{eqconbespprjhlk73} with \eref{eqprbrdemhqjh8ksp31} one obtains
the inequality \eref{eqthmartpreqjhk8kj2} of theorem \ref{SMPrtoshcolatrho}.\\
Combining \eref{eqconbespprjhlk73} with \eref{eqprbrdemhqjh8ksp31} and \eref{eqprbrdemhqjh8ksp30} one obtains
the inequality \eref{eqthmartpreqjhk8kj1} of theorem \ref{SMPrtoshcolatrho} by choosing
$q=\big(\lambda^2(f_1-f_2)t_0 e\big)^{-1}$ ($q>2$ under the condition imposed on $\lambda$).

\subsubsection{Upper bound estimate  \eref{dsjddddjsh99111d1} of
theorem \ref{SMSpCocorkjgptperme1}}\label{dkhdbdjbjdbdb88811}
Theorem \ref{SMPrtoshcolatrho} can be used to give quantitative estimates on any operator as soon as a cell problem is well defined.
Consider $y_t$ is a diffusion on $\R^d$ that may be decomposed  for $t>0$ as
\begin{equation}\label{SMpotohotopoytmx}
y_t=x+\chi(t)+M_t
\end{equation}
where $\chi(t)$ is a uniformly (in $t$) bounded random vector process ($\|\chi\|_\infty \leq C_\chi$) and  $M_t$  is a continuous square integrable $\F_t$ adapted martingale such that $M_0=0$.\\
Assume that for all $l\in\R^d$ with $|l|=1$
there exists a function $f: \R^+ \rightarrow \R^+$ such that   for all $t_2 > t_1 \geq 0$ one has a.s.
\begin{equation}\label{SMpotohotopoytmx000}
 \E[\int \limits_{t_1}^{t_2}d<M.l,M.l>_s | \F_{t_1}] \leq \int \limits_{0}^{t_2-t_1} f(s)ds
\end{equation}
With $f(s) =f_1$ for $s<t_0$ and $f(s)={^tl}Dl<f_1$ for $s \geq t_0$ with $t_0 > 0$.\\
where $D$ is a positive definite symmetric matrix.\\
Assume that the diffusion $y_t$ has symmetric Markovian probability densities $p(t,x,y)$ with respect to the  measure $m(dy)$ such that for all $x,y\in \R^d$ and $t>0$
\begin{equation}\label{SMpotohoptxydiagtd}
p(t,x,y)\leq \frac{C_2}{t^\frac{d}{2}}
\end{equation}
and for $\delta>0$
\begin{equation}\label{SMPrtohoupespacaeqde}
\P_x(|y_t-x|\geq \delta)\leq C_3 e^{-C_4\frac{\delta^2}{t}}
\end{equation}
where $C_2,C_3,C_4$ are constants.\\

\begin{Theorem}\label{SMpotohosupthyd}
Let $y_t$ be the dif{f}usion described above.\\ Then
with $k_1=30 \big(e(f_1-\lambda_{\min}(D))t_0 \big)^\frac{1}{2}/\lambda_{\min}(D)$ and $k_2=30+10d\lambda_{\max}(D)(1+C_4)$
\begin{equation}\label{hgvgestrqestimath0}
 k_1|x-y|<t,\qquad k_2 < \frac{|x-y|}{\sqrt{t}},\qquad |x-y|>4 C_\chi
\end{equation}
one has
\begin{equation}\label{hgvgestrqestimath}
p(t,x,y)\leq \frac{E_1}{t^\frac{d}{2}}\exp\big(-(1-E)\frac{|y-x|^2_{D^{-1}}}{2t}\big)
\end{equation}
with $E_1=C_2(5(\lambda_{\min}(D)C_4)^{-1}+2^d C_3)$
and $E= 3\big((\frac{k_1|x-y|}{t})^2+\frac{\sqrt{t}}{|x-y|}\big)\leq \frac{1}{10}$
\end{Theorem}
\begin{proof}
The estimate on the heat kernel $p(t,x,y)$ will follow from the
chain rule and  decomposing it  the probability of  moving away
from $x$ to "a well chosen set containing $y$
 in  the time $tq$" and its complement. More precisely, writing $e_{y-x}:=(y-x)/|y-x|$ and
$A_{\delta}=\{z\in \R^d\,:\, (z-x).e_{y-x}\geq (1-\delta)|x-y|\}$, using \eref{SMpotohoptxydiagtd} one obtains that for
  $t>0$, $x,y\in \R^d$ and $0<q<1$,
\begin{equation}\label{prtohospptxya}
\begin{split}
p(t,x,y)=&\int_{A_{\delta}}p(tq,x,z)p(t(1-q),z,y)m(dz) \\&+\int_{A^c_{\delta}}p(tq,x,z)p(t(1-q),z,y)m(dz) \\
\leq& \frac{C_2}{t^\frac{d}{2}}\Big[\frac{1}{(1-q)^\frac{d}{2}}\P_x(y_{tq}.e_{y-x}\geq |x-y|(1-\delta))
\\&+\frac{1}{q^\frac{d}{2}}\P_y(|y_{t(1-q)}|\geq \delta
|x-y|)\Big].
\end{split}
\end{equation}
Let's choose $\delta=\exp\big(-|x-y|\big(d D(e_{x-y}) \sqrt{t}\big)^{-1}\big)$ and $q=1-2D(e_{x-y}) C_4 \delta$ (we will use the notation $D(l):={^tl}Dl$).\\
For $|x-y|/\sqrt{t}>\max(dD(e_{x-y})\ln(4D(e_{x-y})C_4),3d D(e_{x-y}))$
(which basically says that the heat kernel is far from its diagonal behavior) one has $\delta < 1/10$ and $1/2<q<1$.
Using the Aronson type estimate \eref{SMPrtohoupespacaeqde} one controls the second term in \eref{prtohospptxya}
\begin{equation}\label{prtohospptxya3}
\P_y(|y_{t(1-q)}|\geq \delta |x-y|) \leq C_3 \exp(-\frac{|x-y|^2}{2D(e_{x-y})t})
\end{equation}
By the properties \eref{SMpotohotopoytmx}, \eref{SMpotohotopoytmx000} and
the corollary \ref{SMPrtoaphoupgoma} one controls the first term in \eref{prtohospptxya}:
 for $r<1$ with  $r=\frac{C_1 \rho}{ qt}$, $\rho=|x-y|(1-\delta)-C_\chi$,
and $C_1=\big(2e(f_1-D(e_{x-y}))t_0 \big)^\frac{1}{2}/D(e_{x-y})$ one has
\begin{equation}\label{prtohospptxya2}
\P_x(y_{tq}.e_{y-x}\geq |x-y|(1-\delta)))\leq e^{\frac{3}{2}r^2}\exp\big(-(1-r^2)\frac{\rho^2}{2 D(e_{x-y})t q} \big)
\end{equation}
Combining \eref{prtohospptxya2}, \eref{prtohospptxya3}, \eref{prtohospptxya} and using
the value of $q$ and $\delta$ given above one easily obtains the estimate
\eref{hgvgestrqestimath} of theorem \ref{SMpotohosupthyd} under the conditions \eref{hgvgestrqestimath0}.
\end{proof}

Now theorem \ref{SMSpCocorkjgptperme1} is a straightforward application of  theorem \ref{SMpotohosupthyd} and a trivial adaptation of the constants appearing in theorem \ref{SMpotohosupthyd}.
Consider $p(t,x,y)$ the heat kernel associated to the Dirichlet form \eref{eqdirforfeeduf1}. Since $p(t,x,y)$ is continuous
 in $L^\infty(\T^d)$ norm with respect to $U$ (we refer to \cite{ChQiHu98} whose result can easily be
adapted to our case)
and $C^\infty(\T^d)$ is dense in $L^\infty(\T^d)$ with respect to
that norm, one can assume $U$ to be smooth and the general result follows by observing that the estimates in
theorem \ref{SMpotohosupthyd} depend only on $\Osc(U)$.\\
By definition $y_t$ has symmetric probability densities with respect
to the measure $m_U$ and the following Aronson type upper
bound is available \cite{Se98}.
\begin{equation}
p(t,x,y) e^{-2U(y)}\leq C
e^{(4+d)\Osc(U)}\frac{1}{t^\frac{d}{2}}\exp\big(-\frac{|x-y|^2}{4t}\big).
\end{equation}
It follows that the conditions   \eref{SMpotohoptxydiagtd} and \eref{SMPrtohoupespacaeqde} are satisfied with constants
 $C_2,C_3,C_4$ depending only on $d$ and $\Osc(U)$. Now write $\chi_l$ the solution of the  associated cell problem: for
$l\in \S^d$, $L_U\chi_l=-l\nabla U$ with $\chi(0)=0$.\\
Using the  theorem 5.4,  chapter 5 of \cite{St1} on elliptic equations with discontinuous coefficients
(see also \cite{St2}), using the periodicity of $\chi$ and observing that $\chi_l(x)=l.x-F_l(x)$
where $F_l$ is harmonic with respect to $L_U$ one easily obtains that
\begin{equation}\label{eqcontrchi}
C_\chi=\|\chi_l\|_\infty \leq C_d \exp\big((3d+2)\Osc(U)\big)
\end{equation}
From Ito formula one has $l.y_t=x+\chi_l(y_t)-\chi_l(x)+\int_0^t (l-\nabla\chi_l)d\omega_s$,
which corresponds to the decomposition given in \eref{SMpotohotopoytmx}.
 The martingale can be written $l.M_t=\int_0^t (l-\nabla\chi_l)d\omega_s$ and
its bracket is equal to $<l.M,l.M>_t=\int_0^t |l-\nabla \chi_l(y_s)|^2 ds$.
It is easy to obtain from the  theorem 3.9 of \cite{GiTr83}  that
\begin{equation}
f_1=\|\nabla \chi_l\|_\infty \leq C_d (1+\|\nabla U\|_\infty) \exp\big((3d+2)\Osc(U)\big)<\infty
\end{equation}
Writing $\phi_l$ the periodic solution of the ergodicity problem $L_U=|l\nabla \chi_l|^2- {^tlD(U)l}$ ($\phi_l(0)=0$) and observing
that $\phi_l=F_l^2-{^tlD(U)l} \psi_l$ where $L_U \psi_l=1$ it is easy to obtain from \eref{eqcontrchi}, the
 theorem 5.4,  chapter 5 of \cite{St1} and the periodicity of $\phi_l$ that
 \begin{equation}
C_\phi=\|\phi_l\|_\infty \leq C_d \exp\big((9d+4)\Osc(U)\big)
\end{equation}
Since, from the Ito formula
\begin{equation}
\E_x[<l.M,l.M>_t]=\E\big[\phi(y_t)-\phi(x)\big]+ t {^tl}D(U)l
\end{equation}
the martingale satisfies the conditions of  theorem \ref{SMpotohosupthyd} with $f_2={^tl}D(U)l$
 and $t_0=C_\phi\big/\big(f_1-\lambda_{\min}(D)\big)$.
Now one can use theorem \ref{SMpotohosupthyd} to obtain a quantitative control on the heat kernel.
It is important to note that all the constants appearing in that theorem only depend and $d$ and $\Osc(U)$
 except may be $k_1=30 \big(e(f_1-\lambda_{\min}(D))t_0 \big)^\frac{1}{2}/\lambda_{\min}(D)$ in which $f_1$ appears.
 This is where the trick operates, indeed
 $(f_1-\lambda_{\min}(D))t_0=C_\phi$ which is a constant depending only on $\Osc(U)$ and $d$.
 Thus in reality all the constants only depends on the dimension and on $\Osc(U)$.
Which proves the upper bound in theorem \ref{SMSpCocorkjgptperme1}.

\subsubsection{Lower bound estimate \eref{dsjddddjsh99111d2} of
 theorem \ref{SMSpCocorkjgptperme1}}
Let $y_t$ the diffusion associated to the Dirichlet form \eref{eqdirforfeeduf1}. As it has been done in  subsection
\ref{dkhdbdjbjdbdb88811} one can prove the estimate \eref{dsjddddjsh99111d2} assuming that $U$ is smooth and the
general case will follow by the continuity of the heat kernel with respect to $U$ in $L^\infty(\T^d)$ norm.\\
First, we will need the following estimate.
\begin{Proposition}\label{SMSpCocorkjgptperme1th1}
For $l \in \S^d$, $\lambda>k_{5,d,\Osc(U)}$ and $k_{6,d,\Osc(U)} \lambda< t$
one has
\begin{equation}
\P[y_t^U.l\geq \lambda] \geq  \frac{1}{4\sqrt{2\pi}}\int_X^\infty e^{-z^2/2}dz
\end{equation}
with $X=\frac{\lambda}{\sqrt{{^tl}D(U)lt}}(1+F)$ and  $F=\frac{k_{7,d,\Osc(U)}}{\lambda}+k_{8,d,\Osc(U)} \sqrt{\frac{\lambda}{t}}\leq \frac{1}{10}$
\end{Proposition}
\begin{proof}
For $l\in \S^d$ let $F_l,\chi_l,\phi_l$ be the functions introduced in \ref{dkhdbdjbjdbdb88811}.
Write $\F_t$ the filtration associated to
Brownian motion appearing in the SDE solved by $y_t$.
 $F_l(y_t)$ is a $(\P,\F_t)$-continuous local martingale vanishing at $0$ such that (Ito calculus)
\begin{equation}\label{eqdecfllobomtandj81}
<F_l,F_l>_t=t D(l)+\phi_l(y_t)+M_t
\end{equation}
with $M_t=-\int_0^t \nabla \phi_l(y_s) d\omega_s$.   Since  $<F_l,F_l>_\infty=\infty$ a.s. by
Dambis, Dubins-Schwarz representation theorem  $B_t=F_l(y_{T_t})$ is a $(\F_{T_t})$-Brownian motion
with $F_l(y_{t})=B_{<F_l,F_l>_t}$ and
\begin{equation}
T_t=\inf \{s\,:\, <F_l,F_l>_s >t\}
\end{equation}
The idea of the proof is then to show that probability of $y_t$ to move away from $0$ behaves like the probability of a BM
of variance $D(l)$ to move away, to achieve this it will be sufficient to show that $M_t$ becomes negligible in front of
$t D(l)$ using the corollary  \ref{SMPrtoaphoupgoma} to control $\P(M_t \geq x)$. More precisely we will use the following lemma.
\begin{Lemma}
For
\begin{equation}\label{eqdecfllobomtandj82}
\lambda>0,\quad \nu>\|\phi_l\|_\infty,\quad \mu>0,\quad \lambda+\|\chi_l\|_\infty+\mu \leq C_2 \mu\sqrt{D(l)t\nu^{-1}}
\end{equation}
 one has
\begin{equation}\label{prtoholocayblmn}
\P[y_t.l\geq \lambda] \geq  1/2 \P[B_{D(l)t} \geq \lambda+\|\chi_l\|_\infty+\mu]-\P[|M_t|\geq \nu-\|\phi_l\|_\infty]
\end{equation}
\end{Lemma}
\begin{proof}
Let $\lambda>0$, from the representation theorem $\P[F_l(y_t)\geq \lambda ]= \P[B_{D(l)t}+E_t\geq \lambda]$ with
$E_t=B_{<F_l,F_l>_t}-B_{D(l)t}$. It follows that for $\mu>0$
\begin{equation}\label{eqprlobomerdui7kjhb1}
\P[F_l(y_t)\geq \lambda]\geq \P[B_{D(l)t} \geq \lambda+\mu]-\P[|E_t|>\mu]
\end{equation}
It follows from \eref{eqdecfllobomtandj81} that for $\nu>0$,
$\P[|E_t|\geq \mu]\leq \P[|\phi(y_t)+M_t|\geq \nu]+\P[\sup_{|z|<\nu}|B_{D(l)t+z}-B_{D(l)t}|\geq \mu]$.
From which one deduces
\begin{equation}\label{eqprlobomerdui7kjhb2}
\P[|E_t|\geq \mu]\leq   \P[|M_t|\geq \nu-\|\phi_l\|_\infty]+2 \P[|B_\nu|> \mu]
\end{equation}
Combining \eref{eqprlobomerdui7kjhb1} and \eref{eqprlobomerdui7kjhb2} one obtains that $\nu>\|\phi_l\|_\infty$
\[
\begin{split}
\P[y_t.l \geq \lambda]\geq& \P[B_{D(l)t} \geq \lambda+\|\chi_l\|_\infty+\mu]
-4\P[B_{D(l)t}\\ \geq &\mu\sqrt{\frac{D(l)t}{\nu}}]-\P[|M_t|\geq \nu-\|\phi_l\|_\infty]
\end{split}\]
Which leads to \eref{prtoholocayblmn} under the last condition in \eref{eqdecfllobomtandj82}.
\end{proof}

Now let us show that
\begin{Lemma}
For  $C_M x<t$ one has
\begin{equation}\label{equgjanashsislsbdg1}
\P(M_t \geq x) \leq  3 \exp\big(-\frac{x^2}{f_2t} \big)
\end{equation}
where $f_2$ and $C_M$ depend only on $d$ and $\Osc(U)$.
\end{Lemma}
\begin{proof}
Write $G(x)=\frac{1}{2}\phi_l^2-\|\phi_l\|_\infty \phi_l$
Since $$L_U G(x)=|\nabla \phi_l|^2- (\|\phi_l\|_\infty -\phi_l)(|\nabla F_l|^2-D(l))$$ one obtains from Ito formula that
\begin{equation}
\begin{split}
\E[<M,M>_t]&\leq 2\|\phi_l\|_\infty \E[\int_0^t|\nabla F_l|^2(y_s)ds+D(l)t]+\|G\|_\infty\\
&\leq 2\|\phi_l\|_\infty (\|\phi_l\|_\infty+2D(l)
t)+2\|\phi_l\|_\infty^2.
\end{split}
\end{equation}
Thus $M_t$ satisfies the conditions of the corollary \ref{SMPrtoaphoupgoma}
with $f_2=4\|\phi_l\|_\infty D(l)$, $f_1=|\nabla \phi_l|^2_\infty$ and $t_0=4 \|\phi_l\|_\infty^2/(f_1-f_2)$, which
leads to \eref{equgjanashsislsbdg1} by observing that $((f_1-f_2)t_0)^{\frac{1}{2}}/f_2$ is upper bounded by a constant
depending only on $\Osc(U)$ and $d$.
\end{proof}

It follows from the equation \eref{prtoholocayblmn} that under the additional conditions,
\begin{equation}\label{prtoholobcond11}
C_M (\nu-\|\phi_l\|_\infty)<t, \quad \text{and}\quad \lambda +\|\chi_l\|_\infty+\mu<C_3 (\nu-\|\phi_l\|_\infty)
\end{equation}
 (where $C_3$ depends only on $d$ and $\Osc(U)$) one has
\begin{equation}\label{prtoholocayblmn2}
\P[y_t.l\geq \lambda] \geq  1/4 \P[B_{D(l)t} \geq
\lambda+\|\chi_l\|_\infty+\mu].
\end{equation}
Choosing $\nu=\|\phi_l\|_\infty+2/C_3\,(\lambda+\|\chi_l\|_\infty+\mu)$ and $$\mu=4(\lambda+\|\chi_l\|_\infty)^{\frac{3}{2}}(C_2\sqrt{D(l)C_3 t})^{-1}$$
 for $\lambda > \|\chi_l\|_\infty$ and $t>C_4(d,\Osc(U))\lambda$ the
conditions \eref{eqdecfllobomtandj82} and \eref{prtoholobcond11} are satisfied and
\begin{equation}
\mu<C_5(d,\Osc(U)) \lambda \sqrt{\frac{\lambda}{t}}\leq \frac{\lambda}{10}
\end{equation}
and it follows from \eref{prtoholocayblmn2} that
\begin{equation}
\P[y_t.l\geq \lambda] \geq  \frac{1}{4} \P[B_{D(l)t} \geq
\lambda(1+C_5 \sqrt{\frac{\lambda}{t}})+\|\chi_l\|_\infty].
\end{equation}
Which proves proposition \ref{SMSpCocorkjgptperme1th1}.
\end{proof}
Now, let $t>0$, $x,y\in \R^d$ and $p(t,x,y)$ be the heat kernel associated to the Dirichlet form \eref{eqdirforfeeduf1}.
Using the chain rule one obtains that for $0<q <1$ and $\delta>0$
\begin{equation}\label{sjshvdhvdhdv771166}
p(t,x,y)\geq C_{d,\Osc(U)}\P_x\big(y_{tq} \in B(y, \delta \sqrt{t})\big) \inf_{z\in B(y, \delta\sqrt{t})} p((1-q) t,z,y)
\end{equation}
It follows by Aronson estimates that
\begin{equation}\label{sjshvdhvdhdv771167}
\begin{split}
p(t,x,y)\geq &C_{d,\Osc(U)}\P_x\big(y_{tq} \in B(y,
\delta\sqrt{t})\big)
\big(t(1-q)\big)^{-\frac{d}{2}}\\&\exp\big(-C_{d,\Osc(U),2}\delta^2
/(1-q)\big).
\end{split}
\end{equation}
Now for $l \in \R^d$ let us define the probability measure $\bar{\P}_x$ as
\begin{equation}
\frac{d \bar{\P}_x}{d \P_x}= \frac{e^{l.y_t}}{\E_x[e^{l.y_t}]}.
\end{equation}
From now we can assume $x:=0$ and we will fix
\begin{equation}\label{sdjdjdbjbd}
l:= D(U)^{-1} y/(qt)
\end{equation}
and assume
\begin{equation}\label{sdjdjdbjbd2}
|l|\leq 1.
\end{equation}
Writing $\bar{\E}_x$ the expectation associated to $\bar{\P}_x$ one has
\begin{equation}\label{sjshvdhvdhdv771168}
\begin{split}
\P_0\big(y_{tq} &\in B(y, \delta \sqrt{t})\big) = \bar{\E}_0\Big[e^{-l.y_{tq}} 1_{y_{tq} \in B(y, \delta\sqrt{t})} \Big]\E_0[e^{l.y_{tq}}]
\\\geq &e^{-y D(U)^{-1} y/(qt)-C_{d,\Osc(U),3}|y|\delta /(qt^{\frac{1}{2}})}
\bar{\P}_0\Big[y_{tq} \in B(y,\delta \sqrt{t})
\Big]\E_0[e^{l.y_{tq}}].
\end{split}
\end{equation}
Now it is trivial to check that the generator of $y_t$ with respect to $\bar{\E}_x$ is
\begin{equation}
\bar{L}=\Delta/2-\nabla U\nabla +l.\nabla
\end{equation}
Let us write $\bar{p}$ the heat kernel associated to that generator, it is trivial to obtain from \eref{sdjdjdbjbd}, \eref{sdjdjdbjbd2} and theorem 1.4 of  \cite{Nor97} that for
$z\in B(y, \delta \sqrt{t})$ one has
\begin{equation}
\bar{p}(tq,0,z)\geq C_{d,\Osc(U),4} (qt)^{-\frac{d}{2}} \exp(-C_{d,\Osc(U),5}\delta^2/q)
\end{equation}
It follows that
\begin{equation}\label{shhsvdhjvdhdv81}
\bar{\P}_0\Big[y_{tq} \in B(y,\delta \sqrt{t})\Big]\geq C_{d,\Osc(U),6}\delta^d q^{-d/2} \exp(-C_{d,\Osc(U),5}\delta^2/q)
\end{equation}
Moreover for $\lambda>0$
\begin{equation}\label{sjshvdhvdhdv771168ff1}
\begin{split}
\E_0[e^{l.y_{tq}}]\geq \P_0[l.y_{tq}\geq \lambda] e^{\lambda}.
\end{split}
\end{equation}
And choosing $\lambda=l D^{-1}(U)l tq$ one easily obtains from proposition \ref{SMSpCocorkjgptperme1th1} that there exists constants $C_1,C_2,C_3$ depending on $d$ and $\Osc(U)$ such that  for $|y|>C_1$ and $C_2|y|<tq$ one has

\begin{equation}\label{sjshvdhvdhdv771168ff2}
\begin{split}
\E_0[e^{l.y_{tq}}]\geq \exp\big(\frac{y D^{-1}(U)y}{2tq}(1-F)\big)
\end{split}
\end{equation}
with
\begin{equation}
F:= C_3 \big(qt/y^2+|y|/(qt)\big)
\end{equation}
Now let us choose
\begin{equation}
q:=1-\exp(-|x-y|t^{-\frac{1}{2}})
\end{equation}
and
\begin{equation}
\delta:=(1-q)^\frac{1}{2}.
\end{equation}
With these values for $q$ and $\delta$ and
combining \eref{sjshvdhvdhdv771168ff2} with \eref{sjshvdhvdhdv771168} and \eref{sjshvdhvdhdv771167} one obtains that for $|y-x|>C_{7,d,\Osc(U)}$ and $C_{8,d,\Osc(U)}|y-x|<t$ one has
\begin{equation}
p_t(x,y)\geq C_{9,d,\Osc(U)} t^{-d/2}\exp(-(1-F_2)|x-y|^2_{D^{-1}(U)}/(2t))
\end{equation}
with
\begin{equation}
F_2:=C_{10,d,\Osc(U)} \big(t/|x-y|^2+|y-x|/(t)\big).
\end{equation}
It is then easy to deduce the lower bound of theorem \ref{SMSpCocorkjgptperme1} by an appropriate shift of the constants.

\subsection{An analytical inequality for sub-harmonic functions}
\subsubsection{The inequality: Theorem \ref{IntthmSMAldTiger4}}
There is no loss of generality by assuming $\Omega$ to be the segment $(0,1)$. We will give a geometrical proof
theorem \ref{IntthmSMAldTiger4} explaining why we expect the existence of an homotopy invariant constant $C_{d,\Omega}$ in
 conjecture \ref{IntSMAldTiger4}.
The theorem \ref{IntthmSMAldTiger4} is proven if the inequality \eref{equuhbusdciujjh1} is true when $\phi$ and $\psi$ are
Green functions $G_{\lambda}(x,z)$ of $-\nabla (\lambda \nabla)$ with Dirichlet condition on $\partial(0,1)$.\\
Let $(x,y)\in (0,1)^2$, $x<y$. Write $\Omega_1=\{z\in \Omega\,:\, \nabla_z G(x,z)\nabla_z G(y,z)<0\}$.
The inequality \eref{equuhbusdciujjh1} is true if
\begin{equation}\label{ujhahcvzshwj8hw60}
-\int_{\Omega_1} \nabla_z G(x,z)\lambda(z)\nabla_z G(y,z)\,dz \leq \int_{\Omega} \nabla_z G(x,z)
\lambda(z)\nabla_z G(y,z)\,dz
\end{equation}
Write $A_x=\{z\in \Omega\,:\, G(x,z)>G(x,y)\}$ and $A_y=\{z\in \Omega\,:\, G(y,z)>G(x,y)\}$. Integrating by parts one
obtains
\begin{equation}\label{ujhahcvzshwj8hw61}
\int_{A_x} \nabla_z G(x,z)\lambda(z)\nabla_z G(y,z)\,dz
=0=\int_{A_y} \nabla_z G(x,z)\lambda(z)\nabla_z G(y,z)\,dz.
\end{equation}
Now the one dimensional specificity shall be used. Since $G(x,z)$ is increasing from $0$ to $x$ and decreasing from
$x$ to $1$, it follows that $\Omega_1=(x,y)$ and $\big(A_x/\Omega_1\big)\cap \big(A_y/\Omega_1\big)=\emptyset$. Combining
this with \eref{ujhahcvzshwj8hw61} one obtains \eref{ujhahcvzshwj8hw60}, which proves the theorem. Let's note that
a simple computation shows that the constant $3$ is sharp.

\subsubsection{Equivalence with the stability of Green functions: Proposition \ref{Propartpreqgrfost}}
Write for $\epsilon \in [0,1]$ $\lambda_\epsilon(x)=e^{U(x)+\epsilon T(x)}$. Write $\psi_\epsilon$ the solution
of $-\nabla (\lambda_\epsilon \nabla \psi_\epsilon)=g$ with Dirichlet condition on $\overline{\Omega}$ and
$g\in C^\infty(\overline{\Omega}),\,g>0$.\\
Assume conjecture \ref{IntSMAldTiger4} to be true, then
proposition \ref{Propartpreqgrfost} is proven if
\begin{equation}\label{hgdsgcguubvz7451}
e^{-C_{d,\Omega} \|T\|_\infty} \leq
\|\frac{\psi_1}{\psi_0}\|_\infty \leq e^{C_{d,\Omega}
\|T\|_\infty}.
\end{equation}
One obtains by differentiation (writing $L_{\lambda_\epsilon}=-\nabla
\lambda_\epsilon \nabla$)
 $L_{\lambda_\epsilon} \partial_\epsilon \psi_\epsilon=-L_{\partial \epsilon \lambda_\epsilon}\psi_\epsilon$.
Which leads by integration by parts to
\begin{equation}
 \partial_\epsilon \psi_\epsilon=-\int_\Omega \nabla_y G_{\lambda_\epsilon}(x,y)
\lambda_\epsilon (y) \nabla_y G_{\lambda_\epsilon}(y,z) T(y)
g_\epsilon(z)\,dy\,dz.
\end{equation}
Using conjecture \ref{IntSMAldTiger4}
\begin{equation}
\begin{split}
 |\partial_\epsilon \psi_\epsilon|
&\leq \|T\|_\infty \int_\Omega |\nabla_y G_{\lambda_\epsilon}(x,y)
\lambda_\epsilon (y) \nabla_y G_{\lambda_\epsilon}(y,z)|  g_\epsilon(z)\,dy\,dz\\
&\leq \|T\|_\infty C_{d,\Omega} \psi_\epsilon.
\end{split}
\end{equation}
And integrating $\partial_\epsilon \ln \psi_\epsilon\leq \|T\|_\infty C_{d,\Omega}$ one obtains the upper bound in
\eref{hgdsgcguubvz7451} (the lower bound being proven in a similar way).\\
Conversely if conjecture \ref{IntSMAldTiger4} is false one can find $\delta>0$ $x,z \in \Omega^2$ and $g$ being
 a smooth approximation of a Dirac around $z$ such that if \\$T(y)=-\operatorname{Sign}\Big(\nabla_y G_{\lambda_\epsilon}(x,y)
\lambda_\epsilon (y) \nabla_y G_{\lambda_\epsilon}(y,z)\Big)$ one has
\begin{equation}
\partial_\epsilon \ln \psi_\epsilon(x)> \|T\|_\infty (1+\delta)
C_{d,\Omega}.
\end{equation}
Which leads to a contradiction with \eref{Newanghgsinh897hj81}.

\subsection{Sub-diffusive behavior}
\subsubsection{Exit times: theorems \ref{IntSMOnedSuHTcoa6}, \ref{IntSMOnedSuHTcoa5} and proposition
 \ref{IntSMOnedSuHTcoahhhbj5}.}
 For $r>1$, write the number of effective scales
\begin{equation}\label{eqexttimesprconttrefsc2}
n_{ef}(r)=\sup\{n \geq 0\,:\, R_n \leq r\}.
\end{equation}
First let us prove that the exit time from $B(0,r)$ is controlled by the homogenization on those first $n_{ef}(r)$ scales:
\begin{Lemma}
\begin{equation}\label{eqexttimesprconttrefsc}
\frac{r^2}{D(V_0^{n_{ef}})} \frac{1}{C_\tau} \leq \E_0[\tau(0,r)]\leq \frac{r^2}{D(V_0^{n_{ef}})} C_\tau
\end{equation}
with $C_\tau=4 e^{6(K_0+K_1\big/ (\rho_{\min}-1))}$.
\end{Lemma}
\begin{proof}
Write $\E^U$, the expectation with respect to the law of probability associated to the generator $1/2\Delta-\nabla U\nabla$. By  theorem \ref{IntthmSMAldTiger4} and  proposition \ref{Propartpreqgrfost} one obtains that
\begin{equation}\label{eqexttimesprconttrefsc5}
e^{-6 \Osc_r(V_{n_{ef}(r)+1}^\infty)}\leq
\E_0[\tau(0,r)]\Big/\E^{V_0^{n_{ef}(r)}}_0[\tau(0,r)] \leq e^{6
\Osc_r(V_{n_{ef}(r)+1}^\infty)}.
\end{equation}
Bounding, $U_{n_{ef}+1}(x)$ by $\Osc(U_n)\leq K_0$ and for $k\geq n_{ef}+2$, $U_k(x)$ by $\|\nabla U_k\|_\infty |x|\leq K_1 |x|/R_k$ one obtains that for $x\in B(0,r)$
\begin{equation}\label{eqexttimesprconttrefsc6}
|V_{n_{ef}(r)+1}^\infty(x)|\leq K_0+K_1\big/ (\rho_{\min}-1).
\end{equation}
Writing $p_{ef}$ corresponds to the maximum number of periods of the scale $n_{ef}$ included in the segment $[0,r]$: $p_{ef}(r)=\sup\{p \geq 1\,:\, p R_{n_{ef}(r)} \leq r\}$; one obtains
\begin{equation}\label{eqexttimesprconttrefsc7}
\begin{split}
\E^{V_0^{n_{ef}(r)}}_0&[\tau(0,p_{ef}(r) R_{n_{ef}(r)})]
 \\&\leq \E^{V_0^{n_{ef}(r)}}_0[\tau(0,r)]\leq \E^{V_0^{n_{ef}(r)}}_0[\tau(0,(p_{ef}(r)+1)
 R_{n_{ef}(r)})].
\end{split}
\end{equation}
Using $\E^{V_0^{n_{ef}(r)}}_0[\tau(0,k R_{n_{ef}(r)})]=(k R_{n_{ef}(r)})^2 \big/D(V_0^{n_{ef}(r)})$, \eref{eqexttimesprconttrefsc5}, \eref{eqexttimesprconttrefsc6} and \eref{eqexttimesprconttrefsc7} one obtains \eref{eqexttimesprconttrefsc}.
\end{proof}
We will need the following mixing lemma
\begin{Lemma}\label{Pr_To_LD_separation2}
Let $(g,f)\in \big(C^1(T^{d}_{1}\big)^2$ and $R \in \N^*$
$$
\Big |\int_{\T^d}g(x)  f(Rx) dx -  \int_{\T^d}g(x) dx \int_{\T^d}f(x) dx\Big | \leq
\|\nabla g\|_\infty/R \int_{\T^d}\big|f\big| dx
$$
\end{Lemma}
\begin{proof}
The proof follows trivially from the following equation
\begin{equation}
\begin{split}
\int_{\T^d}g(x)  f(Rx) dx - & \int_{\T^d}g(x) dx \int_{\T^d}f(x)
dx\\&= \int_{y\in [0,1]^d,x\in \T^d}  f(Rx+y) (g(x+y/R)-g(x)).
\end{split}
\end{equation}
\end{proof}
From lemma \eref{Pr_To_LD_separation2} we will deduce
 a quantitative estimate on the multi-scale effective diffusivities:
\begin{Lemma}
\begin{equation}\label{eqexttimesprconttrefsc3}
\big(\lambda_{\min} e^{-4K_1/\rho_{\min}} \big)^{n}
 \leq D(V^{n-1})\leq \big(\lambda_{\max}e^{4K_1/\rho_{\min}}
 \big)^n.
\end{equation}
\end{Lemma}
\begin{proof}
The proof of \eref{eqexttimesprconttrefsc3} is based on the following functional mixing estimate (obtained from lemma \ref{Pr_To_LD_separation2}):for $U,W\in C^1(\T)$ and $R\in \N^*$ one has
\begin{equation}\label{eqexttimesprconttrefsc8}
\begin{split}
&e^{-\|\nabla W\|_\infty/R} \leq \\&\int_{\T}e^{U(Rx)+W(x)}dx\Big/\Big(\int_{\T}e^{U(x)}dx\int_{\T}e^{W(x)}dx\Big)
 \leq e^{\|\nabla W\|_\infty/R}
\end{split}
\end{equation}
Then by the explicit formula \eref{sdkjcbkjbckjbrkjfbbfffs1} and a straightforward induction on $n$ one obtains that (using \eref{ModsubContUngradUn})
\begin{equation}
\prod_{k=0}^{n-1} \big(e^{4K_1/r_k}\int_{\T}e^{2U_k(x)}dx\int_{\T}e^{-2U_k(x)}dx\big)^{-1}
 \leq D(V^{n-1}_0)
\end{equation}
\begin{equation}
 D(V^{n-1}_0)\leq \prod_{k=0}^{n-1} \big(e^{-4K_1/r_k}
 \int_{\T}e^{2U_k(x)}dx\int_{\T}e^{-2U_k(x)}dx\big)^{-1}.
\end{equation}
Which leads to \eref{eqexttimesprconttrefsc3} by \eref{ModsubDiffCondUniDUn} and \eref{Modsubboundrnrhonmin}.
\end{proof}
Combining \eref{eqexttimesprconttrefsc3} with \eref{eqexttimesprconttrefsc}, \eref{eqexttimesprconttrefsc2} and \eref{Modsubboundrnrhonmin}, one obtains theorem \ref{IntSMOnedSuHTcoa6}.\\
When the medium is self-similar, we will need the following lemma
\begin{Lemma}
\begin{equation}\label{eqexttimesprconttrefsc4}
\lim_{n \rightarrow
\infty}-\frac{1}{n}\ln\big(D(V^{n-1})\big)=\Pr_\rho(2U)+\Pr_\rho(-2U).
\end{equation}
\end{Lemma}
\begin{proof}
The limit \eref{eqexttimesprconttrefsc4} is a direct consequence of the following theorem that is an application of the theory of level-3 large deviations (we refer to \cite{El85} for a sufficient reminder).
\begin{Theorem}\label{PrTothprholco1}
Let $U\in C^\alpha(\T^d)$ (H{\"{o}}lder continuous with exponent $\alpha>0$). Let $R\in \N$, $R\geq 2$. Then
\begin{equation}
\lim_{n\rightarrow \infty}\frac{1}{n}\ln
\int_{\T^d}\exp\big(\sum_{k=0}^{n-1}U(R^k x))dx=\Pr_R(U).
\end{equation}
\end{Theorem}
We have written $\Pr_R$ is the pressure associated to the scaling
shift induced by $R$ on the torus: For $R\in \N/\{0,1\}$ one can
see the torus as a shift space equipped with the transformation
$s_R$
\begin{equation}
\begin{split}
s_R: \T^d &\longrightarrow \T^d \\
x=\sum_{k=1}^{\infty} \frac{x^k}{R^k} &\longrightarrow Rx=\sum_{k=1}^{\infty} \frac{x^{k+1}}{R^k}
\end{split}
\end{equation}
where for each $k$, $x^k$ is a vector in $B=\{0,1,\ldots,R-1\}^d$ and for each $i\in \{1,\ldots,d\}$
$\sum_{k=1}^{\infty} \frac{x^k_i}{R^k}$ is the expression of $x_i$ in base $R$ ($x^k_i\in \{0,\ldots,R-1\}$).\\
Give $B$ with the discrete topology and $B^{\N^*}$ with the product topology. Write $\mu$ the probability measure on $B$ affecting identical weight $1/R^d$ to each element of $B$ and write $\P_\mu$ the associated product measure on $B^{\N^*}$.\\
With respect to the probability space $\big(B^{\N^*},\B(B^{\N^*}),\P_\mu\big)$ the coordinate representation process $x=(x^1,\ldots,x^p,\ldots)$ is a sequence of i.i.d. random variables distributed by $\mu$. When $x$ is seen as an element of the torus $\T^d$ then the probability measure induced by $\mu$ on the torus is the Lebesgue measure.\\
Define the empirical measure $E_n$ associated to the process $x$ by
\begin{equation}
E_n(x,.)=\frac{1}{n} \sum_{k=0}^{n-1}\delta_{s_R^k \operatorname{cycle}(x,n)}
\end{equation}
where $\operatorname{cycle}(x,n)$ is the periodic point in $B^{\N^*}$ obtained by repeating $(x^1,\ldots,x^n)$ periodically. For each $x$, $E_n(x,.)$ is an element of the space $\Me(B^{\N^*})$ of measures on $B^{\N^*}$ and invariant by the shift $s_R$.\\
Then by theorem 9.1.1 of \cite{El85}, $\{Q_n^{(3)}\}$, the $\P_\eta$ distribution on  $\Me(B^{\N^*})$ of the  empirical process $\{E_n\}$ has a large deviation property  with speed $n$ and entropy function $I^{(3)}_\mu$.\\
 We remind that  for $P\in \Me(B^{\N^*})$, $I^{(3)}_\mu(P)=\int_{B^{\N^*}}I^{(2)}_\mu(\tilde{P}) dP$
where $\tilde{P}$ denotes the marginal distribution of $x^1$ associated to $P$ and $I^{(2)}_\mu$ is the relative entropy of $\tilde{P}$ with respect to $\mu$: $I^{(2)}_\mu(\eta)=\int_{B}\ln \frac{d\eta}{d\mu}\;d\mu $.\\
Choosing $U\in C(\T^d)$, H\"{o}lder continuous with exponent $\alpha$, one deduces from the  large deviation property of $\{Q_n^{(3)}\}$ and Varadhan's theorem that
\begin{equation}\label{PrToLetrldVcycE}
\lim_{n\rightarrow \infty}\frac{1}{n}\ln \int_{\T^d}\exp\big(n E_n(x,U)\big)dx=\Pr_R(U)
\end{equation}
where $\Pr_R(U)$ is the pressure of $U$. We remind that
\begin{equation}\label{ProToVaFoPrenent}
\Pr_R(U)=\sup_{P\in \Me_{s_R}(B^{\N^*})} \{\int U\, dP - I^{(3)}_\mu(P)\}
\end{equation}
where $\Me_{s_R}(B^{\N^*})$ is the space of measures on $B^{\N^*}$ invariant by the shift $s_R$.\\
Since  $U$ is H\"{o}lder continuous
\begin{equation}\label{KECJNCJHbubiub87101}
\begin{split}
|n E_n(x,U)-\sum_{k=0}^{n-1}U(R^k x)| &\leq \sum_{k=0}^{n-1}(\frac{C_d}{R^{n-k}})^\alpha \\
&\leq C(d,\alpha) \sum_{k=0}^{\infty}\frac{1}{R^{k \alpha}}\leq
C(d,\alpha,R)<\infty.
\end{split}
\end{equation}
And one obtains theorem \ref{PrTothprholco1} from \eref{KECJNCJHbubiub87101} and \eref{PrToLetrldVcycE}
\end{proof}
Combining \eref{eqexttimesprconttrefsc4} with \eref{eqexttimesprconttrefsc} and \eref{eqexttimesprconttrefsc2}, one obtains theorem \ref{IntSMOnedSuHTcoa5}.\\
Now let us prove proposition \ref{IntSMOnedSuHTcoahhhbj5}.
The basic properties of the pressure can be found in \cite{Ke98} theorem 4.1.10. (note that the definition of the pressure given here differs from the standard one of the topological pressure by a constant that is $d \ln R$, here $\Pr_R(0)=0$). Let's remind that $\Pr_R$ is a convex function on the space of upper semi continuous functions on the torus to $[-\infty,\infty)$ thus $\Pr_R(U)+\Pr_R(-U)\geq 0$.\\ We will now remind the strict convexity of the topological pressure on a well defined equivalence space:\\
To $s_R$ is associated a scaling operator $S_R$ acting on the periodic continuous functions on $\T^d$
\begin{equation}\label{prtoldthscopeq}
\begin{split}
S_R: C(\T^d) &\longrightarrow C(\T^d) \\
\big(x\rightarrow f(x)\big) &\longrightarrow \big(x\rightarrow
f(s_Rx)=f(Rx)\big).
\end{split}
\end{equation}
Write $\I_{S_R}(\T^d)$ the closed subspace of $\C(\T^d)$ generated by the elements $V-S_R^k V$ with $V\in C(\T^d)$ and $k\in \N$. Write $[U]$ the equivalence class of $U$, then by proposition 4.7 of \cite{Rue78} the function
\begin{equation}
\begin{split}
\Pr_R: C(\T^d)/\I_{S_R}(\T^d) &\longrightarrow [-\infty,+\infty) \\
[U] &\longrightarrow \Pr_R(U)
\end{split}
\end{equation}
is well defined on the set of equivalence classes induced by $\I_{S_R}(\T^d)$ on $C(\T^d)$. Moreover
it is strictly convex on the subset
\begin{equation}
\{[U]\in C(\T^d)/\I_{S_R}(\T^d)\,:\,\int_{\T^d}U(x)dx=0\}.
\end{equation}
We will now prove proposition \ref{IntSMOnedSuHTcoahhhbj5}, since for $c\in \R$, $\Pr(U+c)=\Pr(U)+c$, it is sufficient to assume $\int_{\T^d}U(x)\,dx=0$ and show that
\begin{equation}
\Pr_R(2U)+\Pr_R(-2U)=0 \Leftrightarrow \lim_{n\rightarrow \infty}
\frac{1}{n}\|\sum_{k=0}^{n-1}S_{R^k}U\|_{\infty}=0.
\end{equation}
($\Leftarrow$): This implication is easy since
\begin{equation}
0\leq \Pr_R(2U)+\Pr_R(-2U)\leq \lim_{n\rightarrow \infty}
\frac{4}{n}\|\sum_{k=0}^{n-1}S_{R^k}U\|_{\infty}.
\end{equation}
($\Rightarrow$): Assume  $\Pr_R(2U)+\Pr_R(-2U)=0$ then let $\epsilon>0$. Then by the strict convexity of the pressure as described above there exists $W_1,\ldots,W_k \in C(\T^d)$ and $m_1,\ldots,m_k \in \N/\{0,1\}$, $\lambda_1,\ldots,\lambda_k \in R$ such that $W=\sum_{p=1}^k \lambda_p (W_p-S_{R^{m_p}}W_p)$
and
$\|U-W\|_\infty \leq \epsilon$.
Since $\sum_{p=0}^{n-1}S_{R^p}W$ remains bounded it follows that
\begin{equation}
 \lim_{n\rightarrow \infty} \frac{1}{n}\|\sum_{k=0}^{n-1}S_{R^k}U\|_{\infty} \leq
 \epsilon.
\end{equation}
which leads to the proof.

\subsubsection{Mean squared displacement: proposition \ref{propcompcaryt2ujhb1}
 theorem \ref{IntSMOnedsunmsqdthhvh873} }

Let $y_t$ be the solution of \eref{IntModelsubdiffstochdiffequ}. Write
\begin{equation}\label{prtoldtchscopeqhgvaa1}
n_{flu}(t)=\sup\{n\in \N \,:\, R_n^2 \leq t\}
\end{equation}
$n_{flu}$ shall be the number of fluctuating scales that have an influence on the
 mean squared displacement at the time $t$
(the effective scales plus the perturbation scales). Chose the number of perturbation scales to be
\begin{equation}\label{prtoldtchscopeqhgvaa2}
n_{per}=\inf\{n \in \N\,:\, R^2_{n_{flu}-n} e^{14 n K_0} 10^4
\leq t D(V_0^{n_{flu}}) \}.
\end{equation}
We will now prove the following proposition
\begin{Proposition}\label{djdvdhvdhdeu7711}
For
 $\rho_{\min}>10 e^{30 K_1}$ and $t>R_9$, $n_{per}$ is well defined and
\begin{equation}\label{prtoldtchscopeqhgvaa3}
C_1 e^{-8 n_{per}K_0}D(V_0^{n_{flu}})t \leq \E[y_t^2]\leq C_2 e^{8
n_{per}K_0} D(V_0^{n_{flu}})t.
\end{equation}
\end{Proposition}
\begin{proof}
The proof of \eref{prtoldtchscopeqhgvaa3} is based on analytical inequalities that allow to control the
stability of the homogenization process on the smaller scales under the perturbation of larger ones.
More precisely we will first work on an abstract decomposition of $V$ given by \eref{Modsubfracuinfty}
 into effective scales $U$ perturbation scales $P$
and drift scales $T$: $V=U+P+T$ with  $(U,P,T)\in C^\infty(T^1_{R_U}) \times C^\infty(T^1_{R_W})\times C^\infty(\R)$,
 $R_U,R_W \in \N/\{0,1\}$, $R_W/R_U=R_P\in \N^*$ and $W=U+P$ shall correspond to fluctuating scales.\\
Write $\chi^W$ the solution of the cell problem associated to $L_W$
($L_W \chi_W=-\nabla W$, $\chi_W(0)=0$) and $F_W(x)=x-\chi^W(x)$.
Since $F_W$ is harmonic with respect to $L_W=L_V+\nabla T\nabla$ one obtains by Ito formula that
$F_W(y_t)=\int_0^t \nabla F_W (y_s)d\omega_s -\int_0^t \nabla T \nabla F_W(y_s) ds$ from which one obtains that
\begin{equation}\label{eqpryt2lancg851}
\begin{split}
(\frac{1}{2}&-t \|\nabla T\|^2_\infty)\E[\int_0^t |\nabla
F_W(y_s)|^2 ds] \\&\leq \E[F_W^2(y_t)] \leq  2(1+t \|\nabla
T\|^2_\infty)\E[\int_0^t |\nabla F_W(y_s)|^2 ds].
\end{split}
\end{equation}
Write $\chi^P$ the solution of the cell problem associated to $L_P$ and $F^P=x-\chi^P$. We
will show that
\begin{Lemma}
 $F_W=F_P-H^U$ with
\begin{equation}\label{prtolihaass2}
e^{-4\Osc(P)} x^2 \leq  (F^P(x))^2 \leq e^{4\Osc(P)}x^2
\end{equation}
and
\begin{equation}\label{prtokjjlihaass3}
\|H^U\|_\infty \leq 2 (1+4\|\nabla P\|_\infty) e^{2\Osc(P)}
R_W/R_P.
\end{equation}
\end{Lemma}
\begin{proof}
The inequality \eref{prtolihaass2} is a direct consequence of the explicit formula
$F^P(x)=R_W \int_0^x e^{2P(y)}dy\big/\int_0^{R_W} e^{2P(y)}dy$. The inequality \eref{prtokjjlihaass3} follows from
the explicit formula
$$
H^U(x)=R_W \Big(\frac{\int_0^x e^{2P(y)}dy}{\int_0^{R_W} e^{2P(y)}dy}-\frac{\int_0^x
 e^{2(P(y)+U(y))}dy}{\int_0^{R_W} e^{2(P(y)+U(y))}dy}\Big)
$$
noticing that the period of $P$ and $U$ are $R_W$ and $R_W/R_P$ and lemma \ref{Pr_To_LD_separation2}.

\end{proof}
The long time behavior of $\E[\int_0^t |\nabla F_W(y_s)|^2 ds]$ is a perturbation of $D(W)t$ as shown in
 the following lemma
\begin{Lemma}
If
$$R_P> 16 e^{4 \Osc(P)}(\|\nabla P\|_\infty+\|\nabla T\|_\infty)e^{2\|\nabla T\|_\infty/R_P}$$ then for $t>0$
\begin{equation}\label{prtolihaass4}
\begin{split}
- (R_W^2/R_P^2)& (e^{10 \Osc(P)}/R^2_P) 100 e^{4 \|\nabla T\|_\infty/R_P}\\&+(1/6)e^{-4\Osc(P)}D(W) t
 \leq \E[\int_0^t |\nabla F^W|^2(y_s)ds]
\end{split}
\end{equation}
and
\begin{equation}\label{prtolihafeass5}
\begin{split}
\E[\int_0^t |\nabla F^W|^2&(y_s)ds] \leq 6e^{4\Osc(P)} D(W) t\\&+(R_W^2/R_P^2)
 (e^{10 \Osc(P)}/R^2_P) 900e^{4 \|\nabla T\|_\infty/R_P}.
\end{split}
\end{equation}
\end{Lemma}
\begin{proof}
For the proof of \eref{prtolihaass4} and \eref{prtolihafeass5} by scaling one can assume that $R_W=1$ and $R_U=1/R_P$.
Write for $\zeta>0$
\begin{equation}\label{eqlknasdps82n28h1}
\phi_\zeta=2 \int_0^x \frac{e^{2V(y)}}{\int_0^1
e^{2W(y)}dy}\Big[\int_0^y \frac{e^{2(P-T)(z)}} {\int_0^1
e^{2P(z)}dz}dz -\zeta \int_0^y \frac{e^{-2(P+T)(z)}}{\int_0^1
e^{-2P(z)}dz}dz\Big]dy.
\end{equation}
Using lemma \ref{Pr_To_LD_separation2} to separates the scales in \eref{eqlknasdps82n28h1}, it is an easy exercise to obtain that if $R_P> 16 e^{4 \Osc(P)}(\|\nabla P\|_\infty+\|\nabla T\|_\infty)e^{2\|\nabla T\|_\infty/R_P}$ then
\begin{itemize}
\item for $\zeta=6 e^{4 \Osc(P)}$ one has $\sup_{\R}\phi_\zeta \leq 900 \frac{e^{10 \Osc(P)}}
{R^2}e^{4 \|\nabla T\|_\infty/R}$
\item for $\zeta=\frac{e^{-4 \Osc(P)}}{6}$ one has $\inf_{\R} \phi_\zeta \geq
  -100\frac{e^{10\Osc(P)}}{R^2}e^{4\|\nabla T\|_\infty/R }$
\end{itemize}
Observing that $L_V \phi_\zeta = |l-\chi^W_l|^2-\zeta D(W)$
one deduces \eref{prtolihaass4} and \eref{prtolihafeass5} by applying Ito formula.
\end{proof}
Combining \eref{eqpryt2lancg851}, \eref{prtolihaass2}, \eref{prtolihaass4}, \eref{prtolihaass4} and choosing
$U=V_0^{n_{flu}-n_{per}}$, $P=V_{n_{flu}-n_{per}+1}^{n_{flu}}$, $T=V_{n_{flu}+1}^\infty$ ($R_W=R_{n_{flu}}$,
 $R_P=R_{n_{flu}}/R_{n_{flu}-n_{per}}$) and $n_{flu}$ as defined in \eref{prtoldtchscopeqhgvaa1} one obtains that for
$\rho_{min}>C_{K_1,K_0}$
\begin{equation}
D(V_0^{n_{flu}})t e^{-8 n_{per}K_0}/24 - R_{n_{flu}-n_{per}}^2 500
e^{6 n_{per}K_0} \leq \E[y_t^2],
\end{equation}
\begin{equation}
 \E[y_t^2]\leq  (D(V_0^{n_{flu}})t+ R^2_{n_{flu}-n_{per}}) e^{8
 n_{per}K_0}500.
\end{equation}
Which leads to \eref{prtoldtchscopeqhgvaa3} by the choice \eref{prtoldtchscopeqhgvaa2} for $n_{per}$.
\end{proof}
By the uniform control of the ratios \eref{Modsubboundrnrhonmin} one obtains quantitative estimates on
the number of fluctuating and perturbation scales \eref{prtoldtchscopeqhgvaa1} and \eref{prtoldtchscopeqhgvaa2};
combining them with the control \eref{prtoldtchscopeqhgvaa3} and the exponential speed of convergence of the
multi-scale effective diffusivities towards zero \eref{eqexttimesprconttrefsc3}, one obtains
 proposition \ref{propcompcaryt2ujhb1} and  theorem \ref{IntSMOnedsunmsqdthhvh873}.

\subsubsection{Heat kernel tail: theorem  \ref{IntSMOnetrprdeborathi81}}
As it has been done for the mean squared displacement, the proof of  theorem
\ref{IntSMOnetrprdeborathi81} shall follow from an abstract decomposition of the potential $V$. More precisely,
let $R_W \in \N/\{0,1\}$, $(W,T)\in  C^\infty(T^1_{R_W})\times C^\infty(\R)$, ($\|\nabla T\|_\infty<\infty$) and write
$V=W+T$ and $y_t$ the diffusion associated to $L_V$. It has been shown in the proof of proposition \ref{djdvdhvdhdeu7711} that by decomposing
$W$ into $U+P$ where $U$ is of period $R_W/R_P\in \N$ one has for all $t>0$ and all $x \in \R^d$
\begin{equation}\label{prtolihaass4bgenle}
\E_x[\int_0^t |\nabla F^W|^2(y_s)ds]\leq \zeta_2 {D(W)}t+ \frac{R_W^2}{R_P^2} C^\phi_2
\end{equation}
where the constants $C^\phi_2, \zeta_2$ are those given by the
equation \eref{prtolihafeass5}. We will now show  that from the
control \eref{prtolihaass4bgenle} (and $\|\chi^W\|_\infty\leq R_W$
that
 is given by the explicit formula of the solution of the cell problem) one can deduce the following lemma:
\begin{Lemma}\label{ydhdbhdbhcbhcb77112}
For
\begin{equation}\label{prtohafutrprineht4853}
 R_W \leq h/2
\end{equation}
\begin{equation}\label{prtohafutrprineht83}
\|\nabla T\|_\infty 2^3 \big(\zeta_2 D(W)\big)^\frac{1}{2}\leq (h/t)\leq (R_P\big/(R_W \sqrt{C^\phi_2}))
\zeta_2 {D(W)}
\end{equation}
and
\begin{equation}\label{prtohafutrprineht83ineq92}
(R_P\big/(R_W \sqrt{C^\phi_2}))\zeta_2 D(W) e^{-\frac{h^2}{2^{11} \zeta_2 {D(W)}t}} \leq  (h/t)
\end{equation}
one has
\begin{equation}\label{prtohafutrprineht8dss3ineq92}
\P[y_t\geq h] \leq C e^{-\frac{h^2}{2^{9} \zeta_2 {D(W)}t}}.
\end{equation}
\end{Lemma}
\begin{proof}
The proof of \eref{prtohafutrprineht8dss3ineq92} is based on a control of the Laplace transform of $y_t$, more precisely
it is well known that for $\lambda>0$, and $h>0$ one has
$\P[y_t\geq h]\leq \E[e^{\lambda(y_t-h)}]$. Observing that
$y_t=\chi^W(y_t)+\int_0^t \nabla F^W(y_s) d\omega_s-\int_0^t \nabla T.\nabla F^W(y_s) ds$ and using $\|\chi^W\|_\infty\leq R_W$
one deduces by the Cauchy Schwartz inequality  that
\begin{equation}\label{eqcineelhurpaldhhi8k1}
\begin{split}
\P[y_t\geq h] \leq& e^{\lambda( R_W - h)} \E[e^{2 \lambda \int_0^t \nabla F^W(y_s) d\omega_s}]^\frac{1}{2}
\\&\E[e^{2 \sqrt{t}\|\nabla T\|_\infty \lambda \big(\int_0^t |\nabla F^W_l(y_s)|^2
ds\big)^\frac{1}{2}}]^\frac{1}{2}.
\end{split}
\end{equation}
If $X$ is a positive bounded random variable, $\mu'>0$ and $\lambda'>0$ it is easy to show by  integrating by part over
$d\P(X \geq x)$ and using $\P(X \geq x)\leq \E[\exp(\lambda' (X-x))]$  that
\[\E[\exp(\mu' \sqrt{X})]\leq 1+\mu' \exp(\frac{(\mu')^2}{4\lambda'})\sqrt{\frac{\pi}{\lambda'}}\E[\exp(\lambda' X)]\]
Applying this inequality to \eref{eqcineelhurpaldhhi8k1}
 with $X=\int_0^t |\nabla F^W_l(y_s)|^2 ds$, $\lambda'=8\lambda^2$ and
$\mu'=2\lambda\sqrt{t}\|\nabla T\|_\infty$ and observing by Ito formula that
$\E[e^{2 \lambda \int_0^t \nabla F^W_l(y_s) d\omega_s}]
\leq \E[e^{8\lambda^2 \int_0^t |\nabla F^W_l(y_s)|^2 ds}]^\frac{1}{2}$
one obtains
\[\P[y_t\geq h] \leq C e^{\lambda( R_W - h)} e^{\|\nabla T\|^2_\infty t/4
}\E[e^{8\lambda^2 \int_0^t |\nabla F^W(y_s)|^2 ds}]
\]
Now observe that $\int_0^t \nabla F^W_l(y_s) d\omega_s$ satisfies the conditions of theorem
\ref{SMPrtoshcolatrho} with
$f_2=\zeta_2 {D(W)}$, and $t_0 (f_1-f_2)=\frac{R_W^2}{R_P^2} C^\phi_2$. It follows that for
\begin{equation}\label{prtocodeprlarpe129}
8 \lambda^2 \leq (R_P^2)\big/(2eR_W^2C^\phi_2)
\end{equation} one has $$\E[e^{8\lambda^2 \int_0^t |\nabla F^W
(y_s)|^2 ds}]\leq C R_P^4 (e^{8 \lambda^2
 \zeta_2 {D(W)}t})\big/\big(\lambda^4 (C^\phi_2)^2 R_W^4\big)$$
Assuming $ R_W < h/2$ and choosing $\lambda=\frac{h}{32 \zeta_2 {D(W)}t}$ the condition on $\lambda$ in
 \eref{prtocodeprlarpe129} is satisfied under the right inequality in \eref{prtohafutrprineht83} and one obtains
\[\P[y_t\geq h] \leq C e^{-\frac{h^2}
{2^7 \zeta_2 {D(W)}t}} e^{\|\nabla T\|^2_\infty t/4}
\big(R_P^4(\zeta_2 D(W)t)^4\big)\Big/\big(h^4 (C^\phi_2)^2
R_W^4\big).\]
From this the  result
\eref{prtohafutrprineht8dss3ineq92} follows easily by assuming the
 left inequality in \eref{prtohafutrprineht83} (that basically says that the influence of the drift scales
$\|\nabla T\|_\infty$ is small in front of the influence of the fluctuating scales)
 and the condition \eref{prtohafutrprineht83ineq92}.
\end{proof}

Now let's choose $W=V_0^{n_{flu}}$, $P=V_{n_{flu}-n_{per}+1}^{n_{flu}}$, $T=V_{n_{flu}+1}^\infty$ ($R_W=R_{n_{flu}}$,
 $R_P=R_{n_{flu}}/R_{n_{flu}-n_{per}}$) in lemma \ref{ydhdbhdbhcbhcb77112}. For $p\in \N^*$ define the function
\begin{equation}
n_{per}(p)=\inf\{n\in \N\,:\: (R_{p}\big/R_{p-n})e^{-3 nK_0}D(V_0^{p-1})^\frac{1}{2} \geq 2^9 e^{5K_1} \}
\end{equation}
$n_{per}(p)$ corresponds to the number of perturbation scales among $p$ fluctuating scales. We will from now assume that
$\rho_{\min}\geq 2^9 e^{11 K_1}$, which implies that $n_{per}$ is well defined and $1 \leq n_{per}(p) \leq p$. Define
\begin{equation}\label{nflueqdefjhgck87131}
n_{flu}(t/h)=\inf\{n\in \N \,:\, 2^6(K_1\big/R_{n+1})  e^{2 n_{per}(n) K_0} (D(V_0^n))^\frac{1}{2}\leq h/t\}
\end{equation}
 $n_{flu}-n_{per}$ corresponds to the number of fully homogenized scales given $t/h$.
$n_{flu}$ is well defined and greater than $1$  under the following assumption that basically says that homogenization
has started on at least the first scale.
\begin{equation}\label{subdiodprtrcotha2}
(R_2/K_1) e^{2K_0} 2^{-6} \leq t/h
\end{equation}
By the definition of $n_{flu}$ the left inequality in \eref{prtohafutrprineht83} is satisfied.
Using \eref{nflueqdefjhgck87131}, the
right inequality in \eref{prtohafutrprineht83} is implied by the definition of $n_{per}$.
The inequality \eref{prtohafutrprineht4853} is satisfied if $2 R_{n_{flu}} \leq h$; by
 the definition of $n_{flu}$ this is implied by the following inequality that basically says that the heat kernel
behavior is far from its diagonal regime.
\begin{equation}\label{ondsubeqdensprh2tk17}
h^2\big/(D(V_0^{n_{flu}})^\frac{1}{2}t) \geq 2K_1 e^{2K_0}2^6 e^{2 n_{per}(n_{flu}) K_0}
\end{equation}
By the definition of $n_{flu}$ and $n_{per}$, the inequality \eref{prtohafutrprineht83ineq92} is satisfied by the following
inequality that also says that the heat kernel is far from its diagonal regime.
\begin{equation}\label{ondsubeqdensprhweetk17}
2^{14} e^{4 (n_{per}+1)K_0} \ln \Big[R_{n_{flu}+1}\Big]\leq  h^2\big/(D(V_0^{n_{flu}})t)
\end{equation}
With this assumption, it follows by the inequality \eref{prtohafutrprineht8dss3ineq92} that
\begin{equation}
\P[y_t\geq h] \leq C e^{-\frac{h^2}{2^{11} e^{4 n_{per}K_0} D(V_0^{n_{flu}})t}}
\end{equation}
Using the control \eref{eqexttimesprconttrefsc3} on $D(V_0^{n_{flu}})$, and \eref{Modsubboundrnrhonmin} on the ratios
one obtains  theorem  \ref{IntSMOnetrprdeborathi81}. The condition \eref{subdiodprtrcotha2}
 is translated into the first inequality
in \eref{IntSMOnedsutrprdethbishe1} and the conditions \eref{ondsubeqdensprh2tk17}, \eref{ondsubeqdensprhweetk17}
into the second one.

\paragraph*{Acknowledgments}
This research was done at the EPFL in Lausanne.
The author would like to thank G\'{erard} Ben Arous for stimulating discussions; the idea to investigate on the link
between the slow behavior
of a Brownian motion and the presence of
 an infinite number of scales of obstacle comes from his work in geology and the work of M. Barlow and R. Bass
 on the Sierpinski carpet. Thanks
are also due to Hamish Short and to the referee for carefully reading the manuscript and many useful comments.

\end{document}